\documentclass{amsart}
\usepackage{setspace}
\usepackage{a4}
\usepackage{amssymb,amsmath,amsthm,latexsym}
\usepackage{amsfonts}
\usepackage{amsfonts}
\usepackage{graphicx}
\usepackage{textcomp}
\usepackage{cite}
\usepackage{enumerate}
\usepackage[mathscr]{euscript}
\usepackage{mathtools}
\newtheorem{theorem}{Theorem}[section]

\newtheorem{conjecture}[theorem]{Conjecture}
\newtheorem{corollary}[theorem] {Corollary}
\newtheorem{definition}[theorem]{Definition}
\newtheorem{example}[theorem]{Example}
\newtheorem{lemma} [theorem]{Lemma}

\newtheorem{problem}[theorem]{Problem}
\newtheorem{proposition}[theorem]{Proposition}

\title{This is the title}
\usepackage{amssymb}
\usepackage{amssymb}
\usepackage{amssymb}
\usepackage{amssymb}
\usepackage{amsmath}
\usepackage{tikz}
\usepackage{hyperref}
\usepackage{enumerate}
\usepackage{mathtools}
\usepackage{amsmath}
\usepackage{tikz}
\usepackage{amssymb}
\usepackage{amsmath}
\usepackage{tikz}
\usepackage{hyperref}
\usepackage{tikz}
\usepackage{graphicx}
\usepackage{textcomp}
\usetikzlibrary{cd}
\begin{document}
\begin{center}
{\bf\Large{OPERATOR-VALUED p-APPROXIMATE SCHAUDER FRAMES}}\\

\textbf{K. MAHESH KRISHNA} \\
Statistics and Mathematics Unit\\
Indian Statistical Institute, Bangalore Centre\\
Bengaluru,	Karnataka 560 059 India\\
Email: kmaheshak@gmail.com\\
and \\
\textbf{P. SAM JOHNSON} \\
Department of Mathematical and Computational Sciences\\ 
National Institute of Technology Karnataka, Surathkal\\
Mangaluru 575 025, India  \\
Email: sam@nitk.edu.in\\

Date: \today
\end{center}

\hrule
\vspace{0.5cm}
\textbf{Abstract}: We give an operator-algebraic treatment of theory of p-approximate Schuader frames which includes the theory of operator-valued frames by Kaftal, Larson, and Zhang [\textit{Trans. AMS., 2009}], G-frames by Sun [JMAA, 2006],  factorable weak operator-valued frames by Krishna and Johnson [\textit{Annals of FA, 2022}] and p-approximate Schauder frames by Krishna and Johnson [\textit{J. Pseudo-Differ. Oper. Appl, 2021}] as particular cases. We show that a sufficiently rich theory can be developed even for Banach spaces. We achieve this by defining  various concepts and characterizations in  Banach spaces. These include duality, approximate duality, equivalence, orthogonality and stability.



\tableofcontents

\section{Introduction}
Both theoretically and practically successful theory of frames for Hilbert spaces which originated from the seminal work of Duffin and Schaeffer   \cite{CHRISTENSENBOOK, DUFFIN, HEILBOOK, YOUNG, MEYER1} demanded the development of frame theory for Banach spaces. It was in 1991, when Grochenig \cite{GROCHENIG} introduced the notion of frames for Banach spaces, called as Banach frames, in  connection with the theory of atomic decompositions \cite{FEICHTINGERGROCHENIG, FEICHTINGERGROCHENIG2, FEICHTINGERGROCHENIG3}.  In the sequel $ \mathcal{X},  \mathcal{Y} $ always denote Banach  spaces. Dual of $\mathcal{X}$ is denoted by $\mathcal{X}^*$. The identity operator on $ \mathcal{X}$ is denoted by $ I_\mathcal{X}$. Field  real or complex numbers is denoted by $\mathbb{K}$. The space of all bounded linear operators from $ \mathcal{X}$ to $ \mathcal{Y} $ is denoted by $ \mathcal{B}(\mathcal{X}, \mathcal{Y})$. We set $\mathcal{B}(\mathcal{X}) \coloneqq \mathcal{B}(\mathcal{X}, \mathcal{X})$.
\begin{definition}\cite{GROCHENIG}\label{GROCHENIG}
	Let $\mathcal{X}$ be a Banach space,    $\mathcal{X}_d$ be an associated BK-space,  $S:\mathcal{X}_d\to \mathcal{X}$ be a bounded linear operator and  $\{f_n\}_n$ be a collection in   $\mathcal{X}^*$.   The pair $ (\{f_n \}_{n}, S) $ is said to be a \textbf{Banach frame} for $\mathcal{X}$ if 	the following conditions hold.
	\begin{enumerate}[\upshape(i)]
		\item $\{f_n(x)\}_n \in \mathcal{X}_d$, for each $x \in \mathcal{X}$.
		\item There are $a,b>0$ such that 
		\begin{align*}
			a\|x\|\leq \|\{f_n(x)\}_n\|\leq b \|x\|, \quad \forall x \in \mathcal{X}.
		\end{align*}
	Constants $a$ and $b$ are called as Banach frame lower and upper bounds, respectively.
		\item $S\{f_n(x)\}_n=x$, for each $x \in \mathcal{X}$.
	\end{enumerate}
\end{definition}
\begin{definition}\cite{GROCHENIG}
	Let $\mathcal{X}$ be a Banach space,    $\mathcal{X}_d$ be an associated BK-space,  $\{\tau_n\}_n$ be a collection in  $\mathcal{X}$ and  $\{f_n\}_n$ be a collection in  $\mathcal{X}^*.$ The pair $ (\{f_n \}_{n}, \{\tau_n\}_n) $ is said to be an \textbf{atomic decomposition} for $\mathcal{X}$ if 	the following conditions hold.
\begin{enumerate}[\upshape(i)]
	\item $\{f_n(x)\}_n \in \mathcal{X}_d$, for each $x \in \mathcal{X}$.
	\item There are $a,b>0$ such that 
	\begin{align*}
		a\|x\|\leq \|\{f_n(x)\}_n\|\leq b \|x\|, \quad \forall x \in \mathcal{X}.
	\end{align*}
	Constants $a$ and $b$ are called as  lower and upper atomic bounds, respectively.
	\item $\sum_{n=1}^{\infty}f_n(x)\tau_n=x$, for each $x \in \mathcal{X}$.
\end{enumerate}	
\end{definition}
Hilbert space frame theory mainly concentrates  around representing every element of the space using frames as a series but this is missing in the definition of Banach frames. After a decade from the  work of Grochenig, it is Aldroubi, Sun, Tang, Han, Larson, Casazza, Christensen, and Stoeva  \cite{CASAZZAHANLARSON, ALDROUBISUNTANG, CHRISTENSTOEVA, CASAZZACHRISTENSENSTOEVA} analyzed when it is possible to write every element as a series by using a sequence in the space or in the dual space. In \cite{CASAZZACHRISTENSENSTOEVA}, it was particularly analyzed the series representation only by using first two conditions in Definition \ref{GROCHENIG}. Moreover, using the approximation properties of Banach space \cite{CASAZZAAPP}, it was proved by Casazza and Christensen \cite{CASAZZACHRISTENSEN} that there are Banach spaces $\mathcal{X}$ such that there is  no collection $\{f_n\}_n$   in $\mathcal{X}^*$ and no collection in $\{\tau_n\}_n$   in $\mathcal{X}$ satisfying 
\begin{align*}
	x=\sum_{n=1}^\infty
	f_n(x)\tau_n, \quad \forall x \in
	\mathcal{X}.
\end{align*}
Motivated from the dilation theory of frames \cite{KASHINKULIKOVA, CZAJA, HANLARSON}, in 1999, Casazza, Han, and Larson \cite{CASAZZAHANLARSON} introduced three other  notions of frames for Banach spaces called as framing or unconditional framing, projective frame and modeled frame. It is then proved in \cite{CASAZZAHANLARSON}   that all these three notions are equivalent. Around the same time, Terekhin also introduced another notion of frames for Banach spaces, again using dilation viewpoint \cite{TEREKHIN2, TEREKHIN3, TEREKHIN}. 

Important paper of Daubechies and DeVore   \cite{DAUBECHIESDEVORE} led Casazza,  Dilworth, Odell, Schlumprecht, and Zsak \cite{CASAZZADILWORTH} to weaken the unconditional convergence  condition in the definition of framing and they defined Schauder frames for Banach spaces as follows. 
\begin{definition}\cite{CASAZZADILWORTH}\label{SCHAUDERFRAME}
	Let $\{\tau_n\}_n$ be a collection in  $\mathcal{X}$ and 	$\{f_n\}_n$ be a collection in  $\mathcal{X}^*.$ The pair $ (\{f_n \}_{n}, \{\tau_n \}_{n}) $ is said to be a Schauder frame for $\mathcal{X}$ if 
	\begin{align*}
		x=\sum_{n=1}^\infty
		f_n(x)\tau_n, \quad \forall x \in
		\mathcal{X}.
	\end{align*}
\end{definition} 
Definition \ref{SCHAUDERFRAME} has been extended by Thomas, Freeman, Odell,  Schlumprecht, and Zsak in 2014 \cite{FREEMAN, THOMAS}.
\begin{definition}\cite{FREEMAN, THOMAS}\label{ASFDEF}
	Let $\{\tau_n\}_n$ be a collection in  $\mathcal{X}$ and 	$\{f_n\}_n$ be a collection in  $\mathcal{X}^*.$ The pair $ (\{f_n \}_{n}, \{\tau_n \}_{n}) $ is said to be an approximate Schauder frame (ASF) for $\mathcal{X}$ if 
	\begin{align*}
		S_{f, \tau}:\mathcal{X}\ni x \mapsto S_{f, \tau}x\coloneqq \sum_{n=1}^\infty
		f_n(x)\tau_n \in
		\mathcal{X}
	\end{align*}
	is a well-defined bounded linear, invertible operator.
\end{definition}

Again looking back to the theory of frames for Hilbert spaces, following generalizations have been proposed.
\begin{enumerate}[\upshape(i)]
	\item frames for subspaces/fusion frames \cite{CASAZZAFUSION, CASAZZASUBSPACE}.
	\item outer frames \cite{OUTER}.
	\item oblique frames \cite{CHRISTENSENOBLIQUE}.
	\item pseudo frames \cite{LIPSEUDO}.
	\item quasi-projectors \cite{FORNASIERQUASI}.
\end{enumerate}
In 2009, Kaftal, Larson, Zhang \cite{KAFTAL}  introduced the notion of operator-valued frames which unified all the notions mentioned earlier.
\begin{definition}\cite{KAFTAL}\label{KAFTAL}
Let $\mathcal{H}, \mathcal{H}_0$ be Hilbert spaces. 	A collection  $ \{A_n\}_{n} $  in $ \mathcal{B}(\mathcal{H}, \mathcal{H}_0)$ is said to be an \textbf{operator-valued frame}  in $ \mathcal{B}(\mathcal{H}, \mathcal{H}_0)$ if the series 
	\begin{align*}
		S_A\coloneqq \sum_{n=1}^\infty A_n^*A_n
	\end{align*}
	converges in the strong-operator topology on $ \mathcal{B}(\mathcal{H})$ to a  bounded invertible operator.
\end{definition}
 In 2006, Sun \cite{SUN1} introduced G-frames which is equivalent to the notion of operator-valued frames. 
 \begin{definition}\cite{SUN1}
Let   $\mathcal{H}, \mathcal{H}_n $	be Hilbert spaces and let $ A_n \in \mathcal{B}(\mathcal{H}, \mathcal{H}_n)$, for each $n\in \mathbb{N}$. The collection  $ \{A_n\}_{n} $ is said to be a \textbf{G-frame} if there are $a,b>0$ such that 
\begin{align*}
	a\|h\|^2\leq \sum_{n=1}^\infty \|A_nh\|^2\leq b \|h\|^2, \quad \forall h \in \mathcal{H}.
\end{align*}
 \end{definition}
 
 It is noticed by the authors of this paper that there are some other notions of frames even for Hilbert spaces, which will not include in the theory of operator-valued frames. These notions are as follow.
\begin{enumerate}[\upshape(i)]
	\item Framings for Hilbert spaces  \cite{CASAZZAHANLARSON}.
	\item  Schauder frames and approximate  Schauder  frames for Hilbert spaces  \cite{FREEMAN, CASAZZADILWORTH, THOMAS}.   
	\item Atomic decompositions for Hilbert spaces  \cite{CASAZZAHANLARSON}.
	\item cb-frames for Hilbert spaces  \cite{LIUJFA}.
	\item Signed frames for Hilbert spaces  \cite{PENGWALDRON}.  
	\item  Pair  frames for Hilbert space  \cite{FEREYDOONI}.  
	\item Controlled frames for Hilbert spaces  \cite{BALAZSANTOINEGRYBOS}.
	\item K-frames for Hilbert spaces \cite{GAVRUTA}
	\item Semi-frames for Hilbert spaces \cite{ANTOINEBALAZS}.
\end{enumerate}

To include most of  the notions of frames for Hilbert spaces,  the notion of  weak operator-valued frames for Hilbert spaces has been introduced in \cite{KRISHNAJOHNSON5}. Unfortunately, the frame operator in this case may not split and to get a reasonable theory like that of Hilbert space frames, a proper subclass of weak operator-valued frames called as factorable weak operator-valued frames has been introduced  and is defined as follows.
\begin{definition}\cite{KRISHNAJOHNSON5}
Given $n \in \mathbb{N}$, define 
\begin{align*}
	L_n : \mathcal{H}_0 \ni h \mapsto L_nh\coloneqq e_n\otimes h \in  \ell^2(\mathbb{N}) \otimes \mathcal{H}_0,
\end{align*}  where $\{e_n\}_{n} $ is  the standard orthonormal basis for $\ell^2(\mathbb{N})$. 	Let  $ \{A_n\}_{n} $ and  $ \{\Psi_n\}_{n} $ be collections in $ \mathcal{B}(\mathcal{H}, \mathcal{H}_0)$. The pair $( \{A_n\}_{n},  \{\Psi_n\}_{n} )$ is said to be a \textbf{factorable weak operator-valued frame}   in $ \mathcal{B}(\mathcal{H}, \mathcal{H}_0) $   if  the following conditions hold. 
\begin{enumerate}[\upshape(i)]
	\item The series 
	\begin{align*}
	S_{A, \Psi} \coloneqq  \sum_{n=1}^\infty \Psi_n^*A_n
	\end{align*}
	converges in the strong-operator topology on $ \mathcal{B}(\mathcal{H})$ to a  bounded  invertible operator.
	\item The map 
	\begin{align*}
		\theta_A:\mathcal{H} \ni h \mapsto   \theta_A h\coloneqq\sum_{n=1}^\infty L_nA_n h \in \ell^2(\mathbb{N}) \otimes \mathcal{H}_0	
	\end{align*}
is a well-defined bounded linear operator.
\item The map 
\begin{align*}
	\theta_A:\mathcal{H} \ni h \mapsto   \theta_A h\coloneqq\sum_{n=1}^\infty L_n\Psi_n h \in \ell^2(\mathbb{N}) \otimes \mathcal{H}_0	
\end{align*}
is a well-defined bounded linear operator.
\end{enumerate}
\end{definition}
Earlier, we listed various notions of frames for Hilbert spaces. There are corresponding notions for Banach spaces.
\begin{enumerate}[\upshape(i)]
	\item Framings for Banach  spaces  \cite{CASAZZAHANLARSON}.
	\item  Schauder frames and approximate  Schauder  frames for Banach  spaces  \cite{FREEMAN, CASAZZADILWORTH, THOMAS}.   
	\item Atomic decompositions for Banach  spaces  \cite{CASAZZAHANLARSON}.
	\item cb-frames for operator spaces  \cite{LIUJFA}.
\item p-frames for Banach spaces \cite{ALDROUBISUNTANG, CHRISTENSTOEVA}
\item p-approximate Schauder frames for Banach spaces \cite{KRISHNAJOHNSON}.
\item pseudo-Schauder frames for Banach spaces \cite{LIULIUZHENG}.
\item K-frames for Banach spaces \cite{RAMUJOHNSON}.
\item Hilbert-Schauder frames for Banach spaces \cite{LIUHILBERTSCHAUDER}.
\end{enumerate}
Most of the results in frame theory for Hilbert spaces arise by doing analysis in the standard separable Hilbert space $\ell^2(\mathbb{N})$ and switching between a given Hilbert space and $\ell^2(\mathbb{N})$. Therefore a natural space to work  for Banach spaces is the classical sequence space $\ell^p(\mathbb{N})$ for $p\in [1,\infty)$. This is the main motivation for the introduction of p-approximate Schauder frames for Banach spaces. 
\begin{definition}\label{PASFDEF}\cite{KRISHNAJOHNSON}
	An ASF $ (\{f_n \}_{n}, \{\tau_n \}_{n}) $  for $\mathcal{X}$	is said to be a \textbf{p-approximate Schauder frame} (written as p-ASF), $p \in [1, \infty)$ if both the maps 
	\begin{align*}
		 & \theta_f: \mathcal{X}\ni x \mapsto \theta_f x\coloneqq \{f_n(x)\}_n \in \ell^p(\mathbb{N}) \text{ and } \\
		 &\theta_\tau : \ell^p(\mathbb{N}) \ni \{a_n\}_n \mapsto \theta_\tau \{a_n\}_n\coloneqq \sum_{n=1}^\infty a_n\tau_n \in \mathcal{X}
	\end{align*}
	are well-defined bounded linear operators. 
\end{definition}
In this paper, we give operator-algebraic version of both factorable weak operator-valued frames and p-approximate Schauder frames in the setting of Banach spaces. In Section \ref{SECTIONTWO} we define the notion of operator-valued p-approximate Schauder frames and   then show that our notion includes factorable weak operator-valued frames, operator-valued frames/G-frames and p-approximate Schauder frames. We also introduce the notion of operator-valued p-Riesz bases. Through an explicit construction we show that influential Naimark-Han-Larson dilation result can be proved even for the operator-version. We then show that a result based on the result of Holub for Hilbert spaces can be derived for Banach spaces which characterizes operator-valued p-approximate Schauder frames. In Section \ref{DUALITYORTHOGONALITYSECTION} we study two important notions associated with a given frame, namely, duality and orthogonality. We characterize dual frames and use orthogonality to generate new frames. In Section \ref{APPDUALSECTION}, we study the notion of approximate duality and characterize them. Section \ref{SIMILARITYCOMPOSITIONANDTENSORPRODUCT} describes equivalence of frames with an explicit description of operator giving similarity. Section \ref{PERTURBATIONS} derives  perturbation results.  Section \ref{FEICHTINGER} states some conjectures.

\section{Operator-valued p-approximate Schauder frames and their characterizations}\label{SECTIONTWO}
Let $\mathcal{X}$,  $\mathcal{Y}$ be Banach spaces and $p\in [1, \infty)$. By $\ell^p(\mathbb{N})\otimes \mathcal{Y}$, we mean the Banach space obtained by tensoring $\ell^p(\mathbb{N})$  and $ \mathcal{Y}$   w.r.t. any cross norm (see  \cite{DEFANTFLOTET, RYAN, DIESTELGRO} for the tensor product of Banach spaces). In the development of paper \cite{KAFTAL}, certain collection of operators played important role. Here we define corresponding operators for Banach spaces.
 Given $n \in \mathbb{N}$, we define 
\begin{align*}
	L_n : \mathcal{Y} \ni x \mapsto L_nx\coloneqq e_n\otimes x \in  \ell^p(\mathbb{N}) \otimes \mathcal{Y},
\end{align*}  where $\{e_n\}_{n} $ is  the standard Schauder  basis for $\ell^p(\mathbb{N})$. Since the norm is a cross norm, it then follows that $L_n$'s are isometries from  $\mathcal{Y} $ to $ \ell^p(\mathbb{N}) \otimes \mathcal{Y}$. We also  define $\Gamma_n$ as the (unique) tensor product of bounded linear operators
\begin{align*}
	 \ell^p(\mathbb{N})\ni\{a_m\}_m\mapsto a_n \in \mathbb{K}, \quad \mathcal{Y}\ni y \mapsto y \in \mathcal{Y}.
\end{align*} 
Since $\mathbb{K}\otimes \mathcal{Y}$ is isometrically isomorphic to $\mathcal{Y}$, we write codomain of $\Gamma_n$ as $\mathcal{Y}$. We also note that $a_n\otimes y=a_ny$. Hence 
\begin{align*}
	\Gamma_n (\{a_m\}_m \otimes y)=a_ny, \quad \forall \{a_m\}_m \in \ell^p(\mathbb{N}), \forall y \in \mathcal{Y}, \forall n \in \mathbb{N}.
\end{align*}
Therefore, for $  n,m \in \mathbb{N}$ we have  
\begin{align}\label{LEQUATION}
	\Gamma_nL_m =
	\left\{
	\begin{array}{ll}
		I_{\mathcal{Y} } & \mbox{if } n=m \\
		0 & \mbox{if } n\neq m
	\end{array}
	\right.
	~\text{and} \quad 
	\sum\limits_{n=1}^\infty L_n\Gamma_nz=z, \quad \forall z \in \ell^p(\mathbb{N})\otimes \mathcal{Y}.
\end{align}
We now introduce the notion of  operator-valued p-approximate Schauder frames and demonstrate that it includes all known notions of frames for which frame operator factors. 
\begin{definition}
	Let  $ \{A_n\}_{n} $ be collection in $ \mathcal{B}(\mathcal{X}, \mathcal{Y})$ and  $ \{\Psi_n\}_{n} $ be a collection in $ \mathcal{B}(\mathcal{Y}, \mathcal{X})$. The pair $( \{A_n\}_{n},  \{\Psi_n\}_{n} )$ is said to be an \textbf{operator-valued p-approximate Schauder frame} (we write \textbf{operator-valued p-ASF})   in $ \mathcal{B}(\mathcal{X}, \mathcal{Y}) $   if  the following conditions hold.
	\begin{enumerate}[\upshape(i)]
		\item The \textbf{frame operator} 
		\begin{align*}
			S_{A, \Psi}: \mathcal{X}\ni x \mapsto \sum_{n=1}^\infty \Psi_nA_nx \in \mathcal{X}
		\end{align*}
	is a well-defined bounded linear invertible operator. Constants $a,b>0$ satisfying 
	 \begin{align*}
		a\|x\|\leq \left\|\sum_{n=1}^\infty\Psi_nA_nx
	 \right\|\leq b\|x\|, \quad \forall x \in  \mathcal{X}
	\end{align*} 
are called as lower and upper frame bounds, respectively.
		\item The \textbf{analysis operator} 
		\begin{align*}
			\theta_A:\mathcal{X} \ni x \mapsto   \theta_A x\coloneqq\sum_{n=1}^\infty L_nA_n x \in \ell^p(\mathbb{N}) \otimes \mathcal{Y}
		\end{align*}
		is  a well-defined bounded linear operator. Constant $c>0$ satisfying 
		\begin{align*}
			 \left\|\sum_{n=1}^\infty L_nA_nx
			\right\|\leq c\|x\|, \quad \forall x \in  \mathcal{X}
		\end{align*} 
		is called as analysis bound.
		 	\item The \textbf{synthesis operator} 
		\begin{align*}
			\theta_\Psi:\ell^p(\mathbb{N})\otimes \mathcal{Y} \ni z\mapsto	\theta_\Psi z \coloneqq\sum\limits_{n=1}^\infty \Psi_n\Gamma_n z \in \mathcal{X} 
		\end{align*}
is a 	well-defined bounded linear operator. Constant $d>0$ satisfying 
\begin{align*}
	\left\|\sum_{n=1}^\infty \Psi_n\Gamma_nz
	\right\|\leq c\|z\|, \quad \forall z \in  \ell^p(\mathbb{N})\otimes \mathcal{Y} 
\end{align*} 
is called as synthesis bound.
	\end{enumerate}
 If $	S_{A, \Psi}=I_\mathcal{X}$, then the frame is called as a \textbf{Parseval operator-valued p-ASF}.  If we do not demand the invertibility of $	S_{A, \Psi}$, then we say that $( \{A_n\}_{n},  \{\Psi_n\}_{n} )$ is an \textbf{operator-valued p-approximate Bessel sequence} (we write \textbf{operator-valued p-ABS}). 
\end{definition}

\begin{example}
\begin{enumerate}[\upshape(i)]
	\item 	Let $\mathcal{H}, \mathcal{H}_0$ be Hilbert spaces and  $( \{A_n\}_{n},  \{\Psi_n\}_{n} )$ be a factorable   weak-OVF in $ \mathcal{B}(\mathcal{H}, \mathcal{H}_0).$ Define 
	\begin{align*}
	B_n\coloneqq A_n, ~ \Phi_n\coloneqq \Psi_n^*, \quad \forall n \in \mathbb{N}.	
	\end{align*}
Then $(\{B_n\}_{n},  \{\Phi_n\}_{n} )$ is  an operator-valued p-ASF in $ \mathcal{B}(\mathcal{H}, \mathcal{H}_0).$ Note that $\Gamma_n=L_n^*$, $\forall n\in \mathbb{N}$.
	\item Let $\mathcal{X}$ be a Banach space and  $( \{f_n\}_{n},  \{\tau_n\}_{n} )$ be a p-ASF for $\mathcal{X}$. Define 
		\begin{align*}
		A_n\coloneqq f_n, ~ \Psi_n:\mathbb{K} \ni \alpha \mapsto \Psi_n\alpha \coloneqq \alpha \tau_n \in \mathcal{X}, \quad \forall n \in \mathbb{N}.	
	\end{align*}
Then 
\begin{align*}
\sum_{n=1}^{\infty}	\Psi_nA_nx=\sum_{n=1}^{\infty}	\Psi_nf_n(x)=\sum_{n=1}^{\infty}f_n(x)\tau_n, \quad \forall x \in \mathcal{X}.
\end{align*}
Hence $(\{A_n\}_{n},  \{\Psi_n\}_{n} )$ is  an operator-valued p-ASF in $ \mathcal{B}(\mathcal{X}, \mathbb{K}).$
	\item Let $\mathcal{X}$, $\mathcal{Y}$ be Banach spaces and  $U:\mathcal{X} \rightarrow \ell^p(\mathbb{N}) \otimes \mathcal{Y}$ and $ V: \ell^p(\mathbb{N}) \otimes \mathcal{Y}\to \mathcal{X}$ be  bounded linear operators such that $VU$ is bounded invertible.	Define 
	\begin{align*}
		A_n=\Gamma_n U, \quad \Psi_n=VL_n, \quad \forall n \in \mathbb{N}.
	\end{align*}  
Then 
	\begin{align*}
	 \sum_{n=1}^\infty \Psi_nA_nx= \sum_{n=1}^\infty VL_n\Gamma_nUx=V\left(\sum_{n=1}^\infty L_n\Gamma_n\right)Ux=VUx, \quad \forall x \in \mathcal{X}.	
\end{align*}
Hence $(\{A_n\}_{n},  \{\Psi_n\}_{n} )$ is  an operator-valued p-ASF in $ \mathcal{B}(\mathcal{X}, \mathcal{Y}) $.
\item This example is motivated from the Cuntz algebra \cite{CUNTZ}. Let $\mathcal{X}$ be a  Banach space,  $n \in \mathbb{N}$ and $A_1, \dots, A_n, \Psi_1,\dots, \Psi_n: \mathcal{X} \to \mathcal{X}$ be bounded linear operators such that $\sum_{j=1}^{n}\Psi_nA_n=I_\mathcal{X}$. Then $(\{A_j\}_{j=1}^n,  \{\Psi_j\}_{j=1}^n )$ is  an operator-valued p-ASF in $ \mathcal{B}(\mathcal{X}) $.
\end{enumerate}	
\end{example}
It is clear that an operator-valued p-ABS need not be an operator-valued p-ASF. Therefore the next attempt is to see that whether we can convert an operator-valued p-ABS to operator-valued p-ASF. In the case of Hilbert spaces, it was Li and Sun \cite{LISUN}  who showed it is possible to expand every G-Bessel sequence to a G-frame just by adding one element (for a proof in finite dimensional Hilbert spaces, see \cite{CASAZZALEONARD}). Surprisingly, using a result of Casazza and Christensen  \cite{CASAZZACHRISTENSEN}, it is proved in \cite{KRISHNAJOHNSON6}  that for Banach spaces one may not able to convert Bessel sequences to frames. Here is a partial answer for operator-valued p-ASFs. 
\begin{theorem}
	Let $(\{A_n\}_{n},  \{\Psi_n\}_{n} )$ be an operator-valued p-ABS in $ \mathcal{B}(\mathcal{X}, \mathcal{Y}) $. If there are $B \in  \mathcal{B}(\mathcal{X}, \mathcal{Y}) $ and $ \Phi \in \mathcal{B}(\mathcal{Y}, \mathcal{X}) $ such that $I_\mathcal{X}-S_{A,\Psi} =\Phi B$, then 
	\begin{align*}
(\{A_n\}_{n}\cup \{B\},  \{\Psi_n\}_{n}\cup \{\Phi\} )		
	\end{align*}
		is a Parseval  operator-valued p-ASF in $ \mathcal{B}(\mathcal{X}, \mathcal{Y}) $. In particular, if $\mathcal{Y}=\mathcal{X}$ and $I_\mathcal{X}-S_{A,\Psi}$ is a square, then every operator-valued p-ABS can be turned into an operator-valued p-ASF.
\end{theorem}
\begin{proof}
	We only need to observe that 
	\begin{align*}
	\sum_{n=1}^\infty \Psi_nA_nx+\Phi Bx=S_{A,\Psi}x+(I_\mathcal{X}-S_{A,\Psi})x=x	, \quad \forall x \in \mathcal{X}.
	\end{align*}
\end{proof}
In general, it is well known  in Hilbert space frame theory that  it is difficult to invert frame operator and get an expansion of every element of the Hilbert space using frame. Thus one needs an algorithm to approximate every element using frames but without using invertibility of frame operator. First such an algorithm was given by Duffin and Schaeffer  \cite{DUFFIN}  and later two more algorithms were given by Grochenig \cite{GROCHENIGACC}. We now derive an algorithm for operator-valued p-ASFs under a condition.
\begin{proposition}
	Let 	 $( \{A_n\}_{n},  \{\Psi_n\}_{n} )$  be an operator-valued   p-ASF in $ \mathcal{B}(\mathcal{X}, \mathcal{Y})$	 with bounds $a$ and  $b$.  For  $ x \in \mathcal{X}$, define 
	$$ x_0\coloneqq0, ~x_n\coloneqq x_{n-1}+\frac{2}{a+b}S_{A,\Psi}(x-x_{n-1}),\quad \forall n \geq1.$$
	If $\|I_\mathcal{X} -\frac{2}{b+a}S_{A,\Psi}\|\leq \frac{b-a}{b+a}$, 	then 
	$$ \|x_n-x\|\leq \left(\frac{b-a}{b+a}\right)^n\|x\|, \quad\forall n \geq1.$$
	In particular, $x_n\to x $ as $n \to \infty$.
\end{proposition} 
\begin{proof}
	We  observe 
\begin{align*}
	x-x_n&=x-x_{n-1}-\frac{2}{a+b}S_{A, \Psi}(x-x_{n-1})\\
	&=\left(I_\mathcal{X} -\frac{2}{b+a}S_{A, \Psi}\right)(x-x_{n-1})\\
	&=\cdots=\left(I_\mathcal{X} -\frac{2}{b+a}S_{A, \Psi}\right)^nx, \quad\forall x \in \mathcal{X} ,  \forall n \geq 1.
\end{align*}	
Hence 
\begin{align*}
	\|x_n-x\|\leq \left\|I_\mathcal{X} -\frac{2}{b+a}S_{A,\Psi}\right\|^n\|x\|\leq \left(\frac{b-a}{b+a}\right)^n\|x\|, \quad\forall n \geq1.
\end{align*}
\end{proof}
Most important result in the development of theory is the following result.
\begin{theorem}\label{FACTTHEOREM}
	Let  $( \{A_n\}_{n},  \{\Psi_n\}_{n} )$  be an operator-valued   p-ASF in $ \mathcal{B}(\mathcal{X}, \mathcal{Y})$. 
	\begin{enumerate}[\upshape(i)]
		\item The analysis operator $	\theta_A:\mathcal{X} \ni x \mapsto   \theta_A x=\sum_{n=1}^\infty L_nA_n x \in \ell^p(\mathbb{N}) \otimes \mathcal{Y} $	is  injective.
		\item The synthesis operator $	\theta_\Psi:\ell^p(\mathbb{N})\otimes \mathcal{Y} \ni z\mapsto\sum\limits_{n=1}^\infty \Psi_n\Gamma_n z \in \mathcal{X}  $	 surjective.
		 \item  $( \{A_nS_{A,\Psi}^{-1}\}_{n},  \{S_{A,\Psi}^{-1}\Psi_n\}_{n} )$  is  an operator-valued   p-ASF in $ \mathcal{B}(\mathcal{X}, \mathcal{Y})$.
		\item Frame operator factors as  $S_{A,\Psi}=\theta_\Psi\theta_A.$
		\item $ P_{A,\Psi} \coloneqq \theta_A S_{A,\Psi}^{-1} \theta_\Psi:\ell^p(\mathbb{N})\otimes \mathcal{Y} \to \ell^p(\mathbb{N})\otimes \mathcal{Y}$ is a projection onto $ \theta_A(\mathcal{X})$.
	\end{enumerate}	
\end{theorem}
\begin{proof}
	Since $S_{A,\Psi}$ is invertible, we get (i), (ii) and (iii). For (iv), we use Equation (\ref{LEQUATION}) and get 
	\begin{align*}
	\theta_\Psi\theta_Ax&=\left(\sum_{n=1}^\infty \Psi_n\Gamma_n\right)\left(\sum_{k=1}^\infty L_kA_kx\right)=\sum_{n=1}^\infty \Psi_n\left(\sum_{k=1}^\infty \Gamma_nL_kA_kx\right) \\
	&=  \sum_{n=1}^\infty \Psi_nA_nx=S_{A, \Psi}x,  \quad \forall x \in \mathcal{X}.	
	\end{align*}
For (v), 
\begin{align*}
	P_{A,\Psi} ^2=\theta_A S_{A,\Psi}^{-1} \theta_\Psi\theta_A S_{A,\Psi}^{-1} \theta_\Psi=\theta_A S_{A,\Psi}^{-1} I_\mathcal{X}\theta_\Psi=	P_{A,\Psi}.
\end{align*}
\end{proof}
In the next result we are proving that projections give frames to subspaces.
\begin{theorem}\label{FRAMEPROJECTION}
	 Let $\mathcal{Z}$ be a closed complementable subspace of $\mathcal{X}$ and $P:\mathcal{X} \to \mathcal{Z}$ be an onto  projection.
	\begin{enumerate}[\upshape(i)]
		\item If  $( \{A_n\}_{n},  \{\Psi_n\}_{n} )$ is  a Parseval  operator-valued p-ASF in $\mathcal{B}(\mathcal{X}, \mathcal{Y})$, then $( \{A_nP\}_{n},  \{P\Psi_n\}_{n} )$ is  a Parseval  operator-valued p-ASF in $\mathcal{B}(\mathcal{Z}, \mathcal{Y})$.
		\item If  $( \{A_n\}_{n},  \{\Psi_n\}_{n} )$ is  an operator-valued p-ASF in $\mathcal{B}(\mathcal{X}, \mathcal{Y})$, then $( \{A_nP\}_{n},  \{P\Psi_n\}_{n} )$ is  an operator-valued p-ASF in $\mathcal{B}(\mathcal{Z}, \mathcal{Y})$.
	\end{enumerate}
\end{theorem}
\begin{proof}
	\begin{enumerate}[\upshape(i)]
		\item We have to show that the frame operator for $( \{A_nP\}_{n},  \{P\Psi_n\}_{n} )$ is the  identity on $\mathcal{Z}$. Let $z\in \mathcal{Z}$. Then 
		\begin{align*}
		\sum_{n=1}^\infty P\Psi_nA_nPz=	P\left(\sum_{n=1}^\infty \Psi_nA_nz\right)=Pz=z.
		\end{align*}
		\item We observe that the frame operator $S_{AP,P\Psi}$ is $PS_{A,\Psi}P$. Since $S_{A,\Psi}$ is invertible on $\mathcal{X}$, $PS_{A,\Psi}P$ is invertible on $\mathcal{Z}$ with inverse  $PS_{A,\Psi}^{-1}P$.
	\end{enumerate}
\end{proof}
In the theory of Hilbert space frames, the notion which is stronger than frames is that of Riesz basis. We next define the notion of Riesz operator-valued p-ASF. For OVFs, this notion was defined by Kaftal, Larson, and Zhang \cite{KAFTAL} and Sun \cite{SUN1}  and for p-ASFs, this notion was defined in \cite{KRISHNAJOHNSON3}. There is also a notion  of Riesz basis for Banach spaces defined in \cite{ALDROUBISUNTANG}. 
\begin{definition}\label{RIESZOVF}
	\begin{enumerate}[\upshape(i)]
		\item An operator-valued p-ASF $( \{A_n\}_{n},  \{\Psi_n\}_{n} )$  in $\mathcal{B}(\mathcal{X}, \mathcal{Y})$ is said to be  an \textbf{operator-valued p-approximate Riesz basis} (we write operator-valued p-ARB)  if $ P_{A,\Psi}= I_{\ell^p(\mathbb{N})}\otimes I_{\mathcal{Y}}$.
	\item 	A collection  $( \{A_n\}_{n},  \{\Psi_n\}_{n} )$  in $\mathcal{B}(\mathcal{X}, \mathcal{Y})$ is said to be  an \textbf{operator-valued p-approximate Riesz sequence}   if it is an operator-valued p-ARB in  $\mathcal{B}(\mathcal{Z}, \mathcal{Y})$, where 
	\begin{align*}
		\mathcal{Z}\coloneqq \overline{\operatorname{span}}\cup_{n=1}^\infty\Psi_n(\mathcal{Y}).
	\end{align*}
	\end{enumerate}
\end{definition}
Riesz bases play important role in the dilation theory of frames. First dilation result of frames was obtained independently by  Han and Larson  \cite{HANLARSON} and  Kashin and Kulikova \cite{KASHINKULIKOVA} (see \cite{CZAJA} for the history of this theorem).  The result was also known to Daubechies \cite{CASAZZAKOVACEVIC}.  Dilation of framings for Banach spaces were obtained in \cite{CASAZZAHANLARSON}.  This work later led to the notion of general dilation theory of operator-valued measures in \cite{HANLARSONLIULIU, HANLARSONLIULIU2}.   For p-ASFs,  dilation theorem is derived in \cite{KRISHNAJOHNSON3}. For OVFs,     general Naimark-Han-Larson dilation theorem was obtained in \cite{HANLIMENGTANG}. Here  we obtain dilation result for operator-valued p-ASFs. This is basically a kind of converse to Theorem \ref{FRAMEPROJECTION}.

\begin{theorem}(Dilation theorem)\label{OPERATORDILATION}
	Let $( \{A_n\}_{n},  \{\Psi_n\}_{n} )$ be  an operator-valued p-ASF  in $ \mathcal{B}(\mathcal{X}, \mathcal{Y})$. 	Then there exist a Banach space $ \mathcal{X}_1 $ which contains $ \mathcal{X}$ isometrically and  bounded linear operators $B_n:\mathcal{X}_1\rightarrow \mathcal{Y}$, $\Phi_n:\mathcal{Y}\rightarrow \mathcal{X}_1,$ $  \forall n  $ such that $(\{B_n\}_{n} ,\{\Phi_n\}_{n})$ is an operator-valued p-ARB  in $ \mathcal{B}(\mathcal{X}_1, \mathcal{Y})$ and $B_n|_{\mathcal{X}}=  A_n,\forall n \in \mathbb{N}$.
\end{theorem}
\begin{proof}
	Since $ I_{\ell^p(\mathbb{N})\otimes \mathcal{Y}}-P_{A,\Psi}$ is idempotent, $ (I_{\ell^p(\mathbb{N})\otimes \mathcal{Y}}-P_{A,\Psi})(\mathcal{X})$ is a Banach space. Define 
	\begin{align*}
	\mathcal{X}_1\coloneqq \mathcal{X}\oplus 	(I_{\ell^p(\mathbb{N})\otimes \mathcal{Y}}-P_{A,\Psi})(\mathcal{X}).
	\end{align*}
Then the map  $\mathcal{X}\ni x \mapsto x\oplus 0 \in \mathcal{X}_1$  is an isometry. Now define 
	\begin{align*}
	&B_n:\mathcal{X}_1\ni x\oplus y\mapsto A_nx+\Gamma_n( I_{\ell^p(\mathbb{N})\otimes \mathcal{Y}}-P_{A,\Psi}) y \in \mathcal{Y},  \\
	& \Phi_n:\mathcal{Y}\ni y\mapsto \Psi_ny\oplus  (I_{\ell^p(\mathbb{N})\otimes \mathcal{Y}}-P_{A,\Psi})L_n y \in \mathcal{X}_1 , \quad \forall n \in \mathbb{N}.
\end{align*}
Then  we have $B_n|_{\mathcal{X}}=  A_n, \forall n \in \mathbb{N}$. 
We now find 
\begin{align*}
	\theta_B(x\oplus y)&=\sum_{n=1}^\infty L_nA_nx+\sum_{n=1}^\infty L_n\Gamma_n (I_{\ell^p(\mathbb{N})\otimes \mathcal{Y}}-P_{A,\Psi})y\\
	&=\theta_Ax+(I_{\ell^p(\mathbb{N})\otimes \mathcal{Y}}-P_{A,\Psi}) y, \quad \forall  x\oplus y \in \mathcal{X}_1.
\end{align*}
and 
\begin{align*}
	\theta_\Phi z&=\sum_{n=1}^\infty \Psi_n\Gamma_n  z\oplus \sum_{n=1}^\infty (I_{\ell^p(\mathbb{N})\otimes \mathcal{Y}}-P_{A,\Psi})L_n \Gamma_nz\\
	&=\theta_\Psi z\oplus (I_{\ell^p(\mathbb{N})\otimes \mathcal{Y}}-P_{A,\Psi})z, \quad \forall   z \in {\ell^p(\mathbb{N})}\otimes\mathcal{Y}.
\end{align*}
We also notice 
\begin{align*}
(I_{\ell^p(\mathbb{N})\otimes \mathcal{Y}}-P_{A,\Psi})\theta_A=	\theta_A-\theta_AS_{A,\Psi}^{-1}\theta_\Psi \theta_A=\theta_A-\theta_A=0
\end{align*}
and 
\begin{align*}
\theta_\Psi	(I_{\ell^p(\mathbb{N})\otimes \mathcal{Y}}-P_{A,\Psi})=\theta_\Psi-\theta_\Psi\theta_AS_{A,\Psi}^{-1}\theta_\Psi=\theta_\Psi-\theta_\Psi=0.
\end{align*}
Using these expressions we get  the frame operator of $(\{B_n\}_{n} ,\{\Phi_n\}_{n})$ as 
\begin{align*}
	S_{B,\Phi}(x\oplus y)&=\theta_\Phi \theta_B (x\oplus y)=\theta_\Phi (\theta_Ax+(I_{\ell^p(\mathbb{N})\otimes \mathcal{Y}}-P_{A,\Psi}) y)\\
	&=\theta_\Psi (\theta_Ax+(I_{\ell^p(\mathbb{N})\otimes \mathcal{Y}}-P_{A,\Psi}) y)\oplus (I_{\ell^p(\mathbb{N})\otimes \mathcal{Y}}-P_{A,\Psi})(\theta_Ax+(I_{\ell^p(\mathbb{N})\otimes \mathcal{Y}}-P_{A,\Psi}) y)\\
	&=(\theta_\Psi \theta_Ax+\theta_\Psi(I_{\ell^p(\mathbb{N})\otimes \mathcal{Y}}-P_{A,\Psi}) y)\oplus ((I_{\ell^p(\mathbb{N})\otimes \mathcal{Y}}-P_{A,\Psi})\theta_Ax+(I_{\ell^p(\mathbb{N})\otimes \mathcal{Y}}-P_{A,\Psi})^2y)\\
	&=(S_{A,\Psi}x+0)\oplus (0+(I_{\ell^p(\mathbb{N})\otimes \mathcal{Y}}-P_{A,\Psi})y)\\
	&=(S_{A,\Psi}\oplus (I_{\ell^p(\mathbb{N})\otimes \mathcal{Y}}-P_{A,\Psi}))(x\oplus y), \quad \forall  x\oplus y \in \mathcal{X}_1.
\end{align*}
Therefore $S_{B,\Phi}$ invertible with inverse $ S_{A,\Psi}^{-1}\oplus (I_{\ell^p(\mathbb{N})\otimes \mathcal{Y}}-P_{A,\Psi})$. Hence $(\{B_n\}_{n} ,\{\Phi_n\}_{n})$ is an operator-valued p-ASF. Next we need to show that it is Riesz. Let $z\in \ell^p(\mathbb{N})\otimes\mathcal{Y}$. Then 
\begin{align*}
	P_{B,\Phi}z&=\theta_BS_{B,\Phi}^{-1}\theta_\Phi z=\theta_BS_{B,\Phi}^{-1}(\theta_\Psi z\oplus (I_{\ell^p(\mathbb{N})\otimes \mathcal{Y}}-P_{A,\Psi})z)\\
	&=\theta_B(S_{A,\Psi}^{-1}\oplus (I_{\ell^p(\mathbb{N})\otimes \mathcal{Y}}-P_{A,\Psi}))(\theta_\Psi z\oplus (I_{\ell^p(\mathbb{N})\otimes \mathcal{Y}}-P_{A,\Psi})z)\\
	&=\theta_B(S_{A,\Psi}^{-1}\theta_\Psi z\oplus (I_{\ell^p(\mathbb{N})\otimes \mathcal{Y}}-P_{A,\Psi})z)\\
	&=\theta_AS_{A,\Psi}^{-1}\theta_\Psi z+(I_{\ell^p(\mathbb{N})\otimes \mathcal{Y}}-P_{A,\Psi})^2z\\
	&=P_{A,\Psi}z+(I_{\ell^p(\mathbb{N})\otimes \mathcal{Y}}-P_{A,\Psi})z=z.
\end{align*}
\end{proof}

It is possible to characterize  frames for Hilbert spaces  using the standard orthonormal basis $\{e_n\}_n$ for  $\ell^2(\mathbb{N})$  \cite{HOLUB}. This result has been extended in \cite{KRISHNAJOHNSON} which shows that p-ASFs can be characterized using the standard Schauder basis  $\{e_n\}_n$  for $\ell^p(\mathbb{N})$. We further generalize this result for operator-valued p-ASFs.

\begin{theorem}
	A pair  $( \{A_n\}_{n},  \{\Psi_n\}_{n} )$ is an operator-valued p-ASF  in $\mathcal{B}(\mathcal{X}, \mathcal{Y})$ if and only if 
\begin{align*}
	A_n=\Gamma_n U, \quad \Psi_n=VL_n, \quad \forall n \in \mathbb{N},
\end{align*}  
where $U:\mathcal{X} \rightarrow \ell^p(\mathbb{N}) \otimes \mathcal{Y}$ and $ V: \ell^p(\mathbb{N}) \otimes \mathcal{Y}\to \mathcal{X}$ are bounded linear operators such that $VU$ is bounded invertible.	
\end{theorem}
\begin{proof}
	$(\Leftarrow)$ Given  $x\in \mathcal{X}$, Equation (\ref{LEQUATION}) gives 
	\begin{align*}
		S_{A, \Psi}x= \sum_{n=1}^\infty \Psi_nA_nx= \sum_{n=1}^\infty VL_n\Gamma_nUx=V\left(\sum_{n=1}^\infty L_n\Gamma_n\right)Ux=VUx.	
	\end{align*}
	Therefore $S_{A, \Psi}$ is bounded invertible which says $( \{A_n\}_{n},  \{\Psi_n\}_{n} )$ is an operator-valued p-ASF  in $\mathcal{B}(\mathcal{X}, \mathcal{Y})$.
	
		$(\Rightarrow)$ Define $Ux\coloneqq \sum_{n=1}^\infty L_nA_nx$, $\forall x \in \mathcal{X}$,  $Vz\coloneqq \sum_{n=1}^\infty \Psi_n\Gamma_n z$, $\forall z \in \ell^p(\mathbb{N}) \otimes \mathcal{Y}$. Then
		\begin{align*}
			\Gamma_n Ux&=\Gamma_n \left(\sum_{k=1}^\infty L_kA_kx\right)=\sum_{k=1}^\infty \Gamma _nL_kA_kx=A_nx, \quad \forall x \in \mathcal{X}, \forall n\in \mathbb{N}
		\end{align*}
	and 
	\begin{align*}
		VL_nz=\sum_{k=1}^\infty \Psi_k\Gamma_k L_nz=\Psi_nz, \quad \forall z \in \ell^p(\mathbb{N}) \otimes \mathcal{Y}, \forall n\in \mathbb{N}.
	\end{align*}
We are left with showing that $VU$ is invertible. For this we find 
\begin{align*}
VUx=\left(\sum_{n=1}^\infty \Psi_n\Gamma_n\right)\left(\sum_{k=1}^\infty L_kA_kx\right) =  \sum_{n=1}^\infty \Psi_nA_nx=S_{A, \Psi}x,  \quad \forall x \in \mathcal{X}
\end{align*}	
which says $VU$  is invertible.
\end{proof}
\begin{corollary}
		A pair  $( \{A_n\}_{n},  \{\Psi_n\}_{n} )$ is an operator-valued p-ABS  in $\mathcal{B}(\mathcal{X}, \mathcal{Y})$ if and only if 
	\begin{align*}
		A_n=\Gamma_n U, \quad \Psi_n=VL_n, \quad \forall n \in \mathbb{N},
	\end{align*}  
	where $U:\mathcal{X} \rightarrow \ell^p(\mathbb{N}) \otimes \mathcal{Y}$ and $ V: \ell^p(\mathbb{N}) \otimes \mathcal{Y}\to \mathcal{X}$ are bounded linear operators. 
\end{corollary}
\begin{corollary}\label{FIRSTCOROLLARY}
	A pair  $( \{A_n\}_{n},  \{\Psi_n\}_{n} )$ is an operator-valued p-ARB  in $\mathcal{B}(\mathcal{X}, \mathcal{Y})$
if and only if 
\begin{align*}
	A_n=\Gamma_n U, \quad \Psi_n=VL_n, \quad \forall n \in \mathbb{N}
\end{align*}  
where $U:\mathcal{X} \rightarrow \ell^p(\mathbb{N}) \otimes \mathcal{Y}$ and $ V: \ell^p(\mathbb{N}) \otimes \mathcal{Y}\to \mathcal{X}$ are bounded linear operators such that $VU$ is bounded invertible and $ U(VU)^{-1}V =I_{\ell^p(\mathbb{N}) \otimes \mathcal{Y}}$.	
\end{corollary}

\section{Duality and orthogonality}\label{DUALITYORTHOGONALITYSECTION}
Given an operator-valued p-ASF  $( \{A_n\}_{n},  \{\Psi_n\}_{n} )$, consider the frame  $( \{B_n\coloneqq A_nS_{A,\Psi}^{-1}\}_{n},  \{\Phi_n\coloneqq S_{A,\Psi}^{-1}\Psi_n\}_{n} )$ in (iii) of Theorem \ref{FACTTHEOREM}. This frame has the property that $	\theta_\Psi \theta_B=\theta_\Phi\theta_A=I_\mathcal{X}$. In general, there are other frames satisfying this property. This brings us to the notion of dual frames.
\begin{definition}\label{DUALFRAME}
	An operator-valued   p-ASF  $ (\{B_n\}_{n} , \{\Phi_n\}_{n} )$  in $\mathcal{B}(\mathcal{X}, \mathcal{Y})$ is said to be a \textbf{dual}  for  an operator-valued   p-ASF  $  ( \{A_n\}_{n},  \{\Psi_n\}_{n} ) $ in $\mathcal{B}(\mathcal{X}, \mathcal{Y})$ if 
	\begin{align*}
		\theta_\Psi \theta_B=\theta_\Phi\theta_A=I_\mathcal{X}.
	\end{align*}
If both $ (\{B_n\}_{n} , \{\Phi_n\}_{n} )$ and  $  ( \{A_n\}_{n},  \{\Psi_n\}_{n} ) $ are operator-valued p-ABS, then we say that they are dual operator-valued p-ABSs.
\end{definition}  
By using the definition of analysis and synthesis operators we will have the following result.
\begin{proposition}\label{DUALCHAR}
	An operator-valued   p-ASF  $ (\{B_n\}_{n} , \{\Phi_n\}_{n} )$  in  $\mathcal{B}(\mathcal{X}, \mathcal{Y})$ is a dual for an operator-valued   p-ASF  $  ( \{A_n\}_{n},  \{\Psi_n\}_{n} ) $ in  $\mathcal{B}(\mathcal{X}, \mathcal{Y})$ if and only if 
	\begin{align*}
		\sum_{n=1}^\infty\Psi_nB_nx= \sum_{n=1}^\infty\Phi_nA_nx=x, \quad \forall x \in \mathcal{X}.
	\end{align*}
\end{proposition}
The frame $( \{A_nS_{A,\Psi}^{-1}\}_{n},  \{S_{A,\Psi}^{-1}\Psi_n\}_{n} )$  is  called as the canonical dual for $  ( \{A_n\}_{n},  \{\Psi_n\}_{n} )$.  They have the following simple property  whose proof follows from the calculation of its frame operator and we leave the proof.
\begin{theorem}\label{CANONICALDUALFRAMEPROPERTYOPERATORVERSIONWEAK}
	Let $( \{A_n\}_{n},  \{\Psi_n\}_{n} )$ be an  operator-valued   p-ASF in $ \mathcal{B}(\mathcal{X}, \mathcal{Y})$ with frame bounds $ a$ and $ b.$ Then
	\begin{enumerate}[\upshape(i)]
		\item The canonical dual for   the canonical dual   for $( \{A_n\}_{n},  \{\Psi_n\}_{n} )$ is itself.
		\item$ \frac{1}{b}, \frac{1}{a}$ are frame bounds for the canonical dual of $( \{A_n\}_{n},  \{\Psi_n\}_{n} )$.
		\item If $ a, b $ are optimal frame bounds for $( \{A_n\}_{n},  \{\Psi_n\}_{n} )$, then $ \frac{1}{b}, \frac{1}{a}$ are optimal  frame bounds for its canonical dual.
	\end{enumerate} 
\end{theorem} 
Complete description of dual frames  for Hilbert spaces was described  Li \cite{LI}. This description has been generalized in the context of factorable weak OVFs in \cite{KRISHNAJOHNSON5} and p-ASFs in \cite{KRISHNAJOHNSON}. Here we  describe duals which cover all of them.
\begin{lemma}\label{FIRSTLEMMA}
Let $( \{A_n\}_{n},  \{\Psi_n\}_{n} )$ be an operator-valued p-ASF  in $\mathcal{B}(\mathcal{X}, \mathcal{Y})$. Then an operator-valued p-ASF  $ (\{B_n\}_{n} , \{\Phi_n\}_{n} )$  in $\mathcal{B}(\mathcal{X}, \mathcal{Y})$	is a dual for $( \{A_n\}_{n},  \{\Psi_n\}_{n} )$ if and only if 
\begin{align*}
	B_n=\Gamma_n U, \quad \Phi_n=VL_n, \quad \forall n \in \mathbb{N}
\end{align*}
where $U:\mathcal{H} \rightarrow \ell^p(\mathbb{N}) \otimes \mathcal{Y}$  is a bounded right-inverse of $\theta_\Psi $, $V: \ell^p(\mathbb{N}) \otimes \mathcal{Y}\to \mathcal{X}$ is a bounded left-inverse of $\theta_A $ such that $VU$ is bounded invertible.	
\end{lemma}
\begin{proof}
$(\Rightarrow)$ Let $ (\{B_n\}_{n} , \{\Phi_n\}_{n} )$  be a dual operator-valued p-ASF  for  $( \{A_n\}_{n},  \{\Psi_n\}_{n} )$.  Then $\theta_\Psi\theta_B =I_\mathcal{X}= \theta_\Phi\theta_A $. Define $ U\coloneqq\theta_B, V\coloneqq\theta_\Phi.$ Then $U:\mathcal{X} \rightarrow \ell^p(\mathbb{N}) \otimes \mathcal{Y}$ is a   right-inverse of $\theta_\Psi $, $V: \ell^p(\mathbb{N}) \otimes \mathcal{Y}\to \mathcal{X}$ is a left-inverse of $\theta_A $ such that $VU=\theta_\Phi\theta_B=S_{B,\Phi}$ is  invertible. We find
\begin{align*}
	&\Gamma_n Ux=\Gamma_n\left(\sum\limits_{k=1}^\infty L_kB_kx\right)=B_nx, \quad \forall x \in \mathcal{X},\\
	& VL_nz=\sum\limits_{k=1}^\infty\Phi_k\Gamma _kL_nz=\Phi_nz, \quad \forall z \in \ell^p(\mathbb{N}) \otimes \mathcal{Y}, \forall n \in \mathbb{N}.
\end{align*}	
$(\Leftarrow)$ It is clear that $ (\{B_n\}_{n} , \{\Phi_n\}_{n} )$ is an operator-valued p-ASF  in $\mathcal{B}(\mathcal{X}, \mathcal{Y})$ and a calculation gives $\theta_B=U$, $\theta_\Phi=V$. For duality, $\theta_\Phi \theta_A=V \theta_A=I_\mathcal{X} $, $ \theta_\Psi\theta_B=\theta_\Psi U =I_\mathcal{X}$.
\end{proof}

\begin{lemma}\label{SECONDLEMMA}
Let $( \{A_n\}_{n},  \{\Psi_n\}_{n} )$ be an operator-valued p-ASF  in $\mathcal{B}(\mathcal{X}, \mathcal{Y})$. Then
\begin{enumerate}[\upshape(i)]
	\item $R:\mathcal{X} \to \ell^p(\mathbb{N})\otimes\mathcal{Y} $ is a bounded right-inverse of $ \theta_\Psi$ if and only if 
	\begin{align*}
		R=\theta_AS_{A,\Psi}^{-1}+(I_{\ell^p(\mathbb{N})\otimes\mathcal{Y}}-\theta_AS_{A,\Psi}^{-1}\theta_\Psi)U,
	\end{align*}
	where $U :\mathcal{X} \to \ell^p(\mathbb{N})\otimes\mathcal{Y}$ is a bounded linear operator.
	\item $L:\ell^p(\mathbb{N})\otimes\mathcal{Y}\rightarrow \mathcal{X} $ is a bounded left-inverse of $ \theta_A$ if and only if 
	\begin{align*}
		L=S_{A,\Psi}^{-1}\theta_\Psi+V(I_{\ell^p(\mathbb{N})\otimes\mathcal{Y}}-\theta_A S_{A,\Psi}^{-1}\theta_\Psi),
	\end{align*}  
	where $V:\ell^p(\mathbb{N})\otimes\mathcal{Y}\to\mathcal{X}$ is a bounded linear operator.
\end{enumerate}			
\end{lemma}
\begin{proof}
	\begin{enumerate}[\upshape (i)]
		\item $(\Leftarrow)$ Let $U :\mathcal{X} \to \ell^p(\mathbb{N})\otimes\mathcal{Y}$ be a bounded linear operator. Then $ \theta_\Psi(\theta_AS_{A,\Psi}^{-1}+(I_{\ell^p(\mathbb{N})\otimes\mathcal{Y}}-\theta_AS_{A,\Psi}^{-1}\theta_\Psi)U)=I_\mathcal{X}+\theta_\Psi U-\theta_\Psi U=I_\mathcal{X}$. Therefore  $R\coloneqq\theta_AS_{A,\Psi}^{-1}+(I_{\ell^p(\mathbb{N})\otimes\mathcal{Y}}-\theta_AS_{A,\Psi}^{-1}\theta_\Psi)U$ is a bounded right-inverse of $\theta_\Psi$.
		
			$(\Rightarrow)$ Let $R:\mathcal{X} \to \ell^p(\mathbb{N})\otimes\mathcal{Y} $ be a bounded right-inverse of $ \theta_\Psi$. Define $U\coloneqq R$. Then $	\theta_AS_{A,\Psi}^{-1}+(I_{\ell^p(\mathbb{N})\otimes\mathcal{Y}}-\theta_AS_{A,\Psi}^{-1}\theta_\Psi)U=	\theta_AS_{A,\Psi}^{-1}+(I_{\ell^p(\mathbb{N})\otimes\mathcal{Y}}-\theta_AS_{A,\Psi}^{-1}\theta_\Psi)R=\theta_AS_{A,\Psi}^{-1}+R-\theta_AS_{A,\Psi}^{-1}=R$.
	\item $(\Leftarrow)$ Let $V: \ell^p(\mathbb{N})\otimes\mathcal{Y}\rightarrow \mathcal{X}$ be a bounded linear operator. Then $(S_{A,\Psi}^{-1}\theta_\Psi+V(I_{\ell^p(\mathbb{N})\otimes\mathcal{Y}}-\theta_A S_{A,\Psi}^{-1}\theta_\Psi))\theta_A=I_\mathcal{X}+V\theta_A-V\theta_A I_\mathcal{X}=I_\mathcal{X}$. Therefore  $V\coloneqq S_{A,\Psi}^{-1}\theta_\Psi+V(I_{\ell^p(\mathbb{N})\otimes\mathcal{Y}}-\theta_A S_{A,\Psi}^{-1}\theta_\Psi)$ is a bounded left-inverse of $\theta_A$.
				
	$(\Rightarrow)$ Let $ L:\ell^p(\mathbb{N})\otimes\mathcal{Y}\rightarrow \mathcal{X}$ be a bounded left-inverse of $ \theta_A$. Define $V\coloneqq L$. Then $S_{A,\Psi}^{-1}\theta_\Psi+V(I_{\ell^p(\mathbb{N})\otimes\mathcal{Y}}-\theta_A S_{A,\Psi}^{-1}\theta_\Psi) =S_{A,\Psi}^{-1}\theta_\Psi+L(I_{\ell^p(\mathbb{N})\otimes\mathcal{Y}}-\theta_A S_{A,\Psi}^{-1}\theta_\Psi)=S_{A,\Psi}^{-1}\theta_\Psi+L-I_{\mathcal{X}}S_{A,\Psi}^{-1}\theta_\Psi= L$. 	
	\end{enumerate}
\end{proof}

\begin{theorem}\label{ALLDUALPOVF}
Let $( \{A_n\}_{n},  \{\Psi_n\}_{n} )$ be an operator-valued p-ASF  in $\mathcal{B}(\mathcal{X}, \mathcal{Y})$. Then an operator-valued p-ASF  $ (\{B_n\}_{n} , \{\Phi_n\}_{n} )$   in $\mathcal{B}(\mathcal{X}, \mathcal{Y})$ is a dual  for $( \{A_n\}_{n},  \{\Psi_n\}_{n} )$ if and only if
\begin{align*}
	&B_n=A_nS_{A,\Psi}^{-1}+\Gamma_nU-A_nS_{A,\Psi}^{-1}\theta_\Psi U,\\
	&\Phi_n=S_{A,\Psi}^{-1}\Psi_n+VL_n-V\theta_AS_{A,\Psi}^{-1}\Psi_n, \quad \forall n \in \mathbb{N}
\end{align*}
such that the operator 
\begin{align*}
S_{A, \Psi}^{-1}+VU-V\theta_AS_{A, \Psi}^{-1}\theta_\Psi U
\end{align*}
is bounded invertible, where   $U :\mathcal{X} \to \ell^p(\mathbb{N})\otimes\mathcal{Y}$ and $V:\ell^p(\mathbb{N})\otimes\mathcal{Y}\to\mathcal{X}$ are bounded linear operators.	
\end{theorem}
\begin{proof}
	Lemmas \ref{FIRSTLEMMA} and  \ref{SECONDLEMMA}
give the characterization of dual operator-valued p-ASF for $( \{A_n\}_{n},  \{\Psi_n\}_{n} )$ as 
\begin{align*}
	&B_n=\Gamma_n(\theta_AS_{A,\Psi}^{-1}+(I_{\ell^p(\mathbb{N})\otimes\mathcal{Y}}-\theta_AS_{A,\Psi}^{-1}\theta_\Psi)U)=A_nS_{A,\Psi}^{-1}+\Gamma_nU-A_nS_{A,\Psi}^{-1}\theta_\Psi U,\\
	&\Phi_n=(S_{A,\Psi}^{-1}\theta_\Psi+V(I_{\ell^p(\mathbb{N})\otimes\mathcal{Y}}-\theta_A S_{A,\Psi}^{-1}\theta_\Psi))L_n=S_{A,\Psi}^{-1}\Psi_n+VL_n-V\theta_AS_{A,\Psi}^{-1}\Psi_n, \quad \forall n \in \mathbb{N}
\end{align*}	
	such that the operator 
	\begin{align*}
		(S_{A,\Psi}^{-1}\theta_\Psi+V(I_{\ell^p(\mathbb{N})\otimes\mathcal{Y}}-\theta_A S_{A,\Psi}^{-1}\theta_\Psi))(\theta_AS_{A,\Psi}^{-1}+(I_{\ell^p(\mathbb{N})\otimes\mathcal{Y}}-\theta_AS_{A,\Psi}^{-1}\theta_\Psi)U)
	\end{align*}
	is bounded invertible, where $U :\mathcal{X} \to \ell^p(\mathbb{N})\otimes\mathcal{Y}$ and $V:\ell^p(\mathbb{N})\otimes\mathcal{Y}\to\mathcal{X}$ are bounded linear operators. We make an expansion and get 
	\begin{align*}
		&(S_{A,\Psi}^{-1}\theta_\Psi+V(I_{\ell^p(\mathbb{N})\otimes\mathcal{Y}}-\theta_A S_{A,\Psi}^{-1}\theta_\Psi))(\theta_AS_{A,\Psi}^{-1}+(I_{\ell^p(\mathbb{N})\otimes\mathcal{Y}}-\theta_AS_{A,\Psi}^{-1}\theta_\Psi)U)\\
		&\quad =S_{A, \Psi}^{-1}+VU-V\theta_AS_{A, \Psi}^{-1}\theta_\Psi U.	
	\end{align*}
\end{proof}	
Notion which is opposite to the notion of duality is the notion of orthogonality of frames 	which  is very useful in the construction of new frames from two frames.  Balan  \cite{BALANTHESIS} introduced this notion and  later studied by Han and Larson \cite{HANLARSON}. A detailed study in the finite dimensional case can be found in  \cite{KORNELSON}. Orthogonality will be defined as follows based in Definition \ref{DUALFRAME}.
\begin{definition}
	An operator-valued   p-ASF  $ (\{B_n\}_{n} , \{\Phi_n\}_{n} )$  in $\mathcal{B}(\mathcal{X}, \mathcal{Y})$ is said to be \textbf{orthogonal} to an operator-valued   p-ASF  $  ( \{A_n\}_{n},  \{\Psi_n\}_{n} ) $ in $\mathcal{B}(\mathcal{X}, \mathcal{Y})$ if 
	\begin{align*}
	\theta_\Psi \theta_B=\theta_\Phi\theta_A=0.
	\end{align*}
\end{definition}  
Similar to Proposition \ref{DUALCHAR}, we have the following result for orthogonality.
\begin{proposition}
	An operator-valued   p-ASF  $ (\{B_n\}_{n} , \{\Phi_n\}_{n} )$  in  $\mathcal{B}(\mathcal{X}, \mathcal{Y})$ is  orthogonal to an operator-valued   p-ASF  $  ( \{A_n\}_{n},  \{\Psi_n\}_{n} ) $ in  $\mathcal{B}(\mathcal{X}, \mathcal{Y})$ if and only if 
	\begin{align*}
		\sum_{n=1}^\infty\Psi_nB_nx= \sum_{n=1}^\infty\Phi_nA_nx=0, \quad \forall x \in \mathcal{X}.
	\end{align*}
\end{proposition}
As mentioned earlier, we can get new frames using orthogonal frames. Following two results will describe them.
\begin{proposition}
	Let $  ( \{A_n\}_{n},  \{\Psi_n\}_{n} ) $ and $ (\{B_n\}_{n} , \{\Phi_n\}_{n} )$ be  two Parseval  operator-valued p-ASF in   $\mathcal{B}(\mathcal{X}, \mathcal{Y})$ which are  orthogonal. If $C,D,E,F \in \mathcal{B}(\mathcal{X})$ are such that $ EC+FD=I_\mathcal{X}$, then  
	\begin{align*}
	 (\{A_nC+B_nD\}_{n}, \{E\Psi_n+F\Phi_n\}_{n})
	\end{align*}
	 is a  Parseval operator-valued p-ASF in  $\mathcal{B}(\mathcal{X}, \mathcal{Y})$. In particular,  if scalars $ c,d,e,f$ satisfy $ce+df =1$, then $ (\{cA_n+dB_n\}_{n}, \{e\Psi_n+f\Phi_n\}_{n}) $ is   a Parseval  operator-valued p-ASF.
\end{proposition} 
\begin{proof}
For each $x \in \mathcal{X}$, 
	\begin{align*}
	S_{AC+BD,E\Psi +F\Phi } x&=\sum_{n=1}^\infty(E\Psi_n+F\Phi_n)(A_nC+B_nD)x\\
	&=ES_{A,\Psi}Cx+E\left(\sum_{n=1}^\infty\Psi_nB_n\right)Dx+F\left(\sum_{n=1}^\infty\Phi_nA_n\right)Cx+FS_{B,\Phi}Dx\\
	&=ECx+0+0+FDx=x.
	 \end{align*}
\end{proof}

\begin{proposition}
	If $  ( \{A_n\}_{n},  \{\Psi_n\}_{n} ) $  and $ (\{B_n\}_{n} , \{\Phi_n\}_{n} )$ are   orthogonal operator-valued p-ASFs in $ \mathcal{B}(\mathcal{X}, \mathcal{Y})$, then  $(\{A_n\oplus B_n\}_{n},\{\Psi_n\oplus \Phi_n\}_{n})$ is an operator-valued  p-ASF in $ \mathcal{B}(\mathcal{X}\oplus \mathcal{X}, \mathcal{Y}).$    Further, if both $  ( \{A_n\}_{n},  \{\Psi_n\}_{n} ) $  and $ (\{B_n\}_{n} , \{\Phi_n\}_{n} )$ are  Parseval, then $(\{A_n\oplus B_n\}_{n},\{\Psi_n\oplus \Phi_n\}_{n})$ is Parseval.
\end{proposition}
\begin{proof}
	Let $ x \oplus x_1 \in \mathcal{X}\oplus \mathcal{X}$. Then 
	\begin{align*}
	S_{A\oplus B, \Psi\oplus \Phi}(x\oplus x_1)&=\sum_{n=1}^\infty(\Psi_n\oplus \Phi_n)(A_n\oplus B_n)( x\oplus  x_1)=\sum_{n=1}^\infty(\Psi_n\oplus \Phi_n)(A_nx+ B_nx_1)\\
	&=\sum_{n=1}^\infty(\Psi_n(A_nx+B_n x_1)\oplus \Phi_n(A_nx+B_n x_1))
	\\
	&=\left(\sum_{n=1}^\infty\Psi_nA_nx+\sum_{n=1}^\infty\Psi_nB_n x_1\right)\oplus \left(\sum_{n=1}^\infty\Phi_nA_nx+\sum_{n=1}^\infty\Phi_nB_n x_1\right)\\
	&=(S_{A,\Psi}x+0)\oplus(0+S_{B,\Phi}x_1) =(S_{A,\Psi}\oplus S_{B,\Phi})( x\oplus  x_1).
	\end{align*}	
\end{proof}  
In the case of Hilbert spaces, it is known that orthogonal operator-valued frames have a common dual operator-valued frame. However, for Banach spaces, we have the following weaker result.
\begin{proposition}
	Two orthogonal operator-valued p-ASFs have a common dual operator-valued p-ABS.
\end{proposition}
\begin{proof}
Let $  ( \{A_n\}_{n},  \{\Psi_n\}_{n} ) $  and $ (\{B_n\}_{n} , \{\Phi_n\}_{n} )$ be  orthogonal operator-valued p-ASFs in $ \mathcal{B}(\mathcal{X}, \mathcal{Y})$. Define 
\begin{align*}
	C_n\coloneqq A_nS_{A,\Psi}^{-1}+B_nS_{B,\Phi}^{-1}, \quad  \Xi_n\coloneqq S_{A,\Psi}^{-1}\Psi_n+S_{B,\Phi}^{-1}\Phi_n\quad \forall n \in \mathbb{N}.
\end{align*}	
Clearly $ (\{C_n\}_{n} , \{\Xi_n\}_{n} )$ is an  operator-valued p-ABS in $ \mathcal{B}(\mathcal{X}, \mathcal{Y})$. We need to show that this is a common dual for $  ( \{A_n\}_{n},  \{\Psi_n\}_{n} ) $  and $ (\{B_n\}_{n} , \{\Phi_n\}_{n} )$. Let $x\in \mathcal{X}$. Then 

\begin{align*}
	&\sum_{n=1}^{\infty}\Xi_nA_nx=	\sum_{n=1}^{\infty}(S_{A,\Psi}^{-1}\Psi_n+S_{B,\Phi}^{-1}\Phi_n)A_nx=x+0=x,\\
	&\sum_{n=1}^{\infty}\Psi_nC_nx=\sum_{n=1}^{\infty}\Psi_n(A_nS_{A,\Psi}^{-1}+B_nS_{B,\Phi}^{-1})x=x+0=x,\\
	&\sum_{n=1}^{\infty}\Xi_nB_nx=\sum_{n=1}^{\infty}(S_{A,\Psi}^{-1}\Psi_n+S_{B,\Phi}^{-1}\Phi_n)B_nx=0+x=x,\\
	&\sum_{n=1}^{\infty}\Phi_nC_nx=\sum_{n=1}^{\infty}\Phi_n(A_nS_{A,\Psi}^{-1}+B_nS_{B,\Phi}^{-1})x=0+x=x.
\end{align*}
\end{proof}
We end this section by showing that there is no need of any conditions on two operator-valued p-ASFs to make the tensor product as an operator-valued p-ASFs. For frames for Hilbert spaces, this result was proved by Feichtinger and Grochenig \cite{FEICHTINGERGROCHENIGTENSOR}.
\begin{proposition}
	If $  ( \{A_n\}_{n},  \{\Psi_n\}_{n} ) $  and $ (\{B_n\}_{n} , \{\Phi_n\}_{n} )$ are    operator-valued p-ASFs in $ \mathcal{B}(\mathcal{X}, \mathcal{Y})$, then  $(\{A_n\otimes B_m\}_{n,m},\{\Psi_n\otimes \Phi_m\}_{n,m})$ is an operator-valued  p-ASF in $ \mathcal{B}(\mathcal{X}\otimes \mathcal{X}, \mathcal{Y}).$    Further, if both $  ( \{A_n\}_{n},  \{\Psi_n\}_{n} ) $  and $ (\{B_n\}_{n} , \{\Phi_n\}_{n} )$ are  Parseval, then $(\{A_n\otimes B_m\}_{n,m},\{\Psi_n\otimes \Phi_m\}_{n,m})$ is Parseval.	
\end{proposition}
\begin{proof}
	We need to show that the frame operator   for the tensor product is invertible. We first find it at elementary tensors as follows.  
	\begin{align*}
		S_{A\otimes B, \Psi \otimes \Phi}(x\otimes x_1)&=\sum_{n,m=1}^{\infty}(\Psi_n \otimes \Phi_m)(A_n\otimes B_m)(x\otimes x_1)\\
		&=\sum_{n,m=1}^{\infty}(\Psi_n A_n x\otimes \Phi_mB_m x_1)=\left(\sum_{n=1}^{\infty}\Psi_n A_n x\right)\otimes \left(\sum_{m=1}^{\infty}\Phi_m B_m x_1\right)\\
		&=S_{A,\Psi}x\otimes S_{B,\Phi}x_1=(S_{A,\Psi}\otimes S_{B,\Phi})(x\otimes x_1), \quad \forall x,x_1\in \mathcal{X}.
	\end{align*}
We therefore have $	S_{A\otimes B, \Psi \otimes \Phi}=S_{A,\Psi}\otimes S_{B,\Phi}$ which completes the proof.
\end{proof}

\section{Approximate duality}\label{APPDUALSECTION}
Even though all duals   for an operator-valued p-ASF   are  characterized completely in the previous section,  even in the case of Hilbert spaces, it is known that given a collection, it is difficult to verify whether it is a dual of a frame. This led  Christensen  and Laugesen  \cite{CHRISTENSENLAUGESEN} to introduce the notion of approximate duals (also see \cite{LIYAN}). This notion has then been introduced and characterized for Banach spaces in \cite{KRISHNAJOHNSON4}. Here we give operator-valued version of those results for Banach spaces. 
\begin{definition}\label{APPROXIMATELYDUALPABS}
	Let  $( \{A_n\}_{n},  \{\Psi_n\}_{n} )$ and  $( \{B_n\}_{n},  \{\Phi_n\}_{n} )$ be operator-valued  p-ABSs in $ \mathcal{B}(\mathcal{X}, \mathcal{Y})$. We say that they are \textbf{approximately dual p-ABSs} if 
	\begin{align*}
		\|I_\mathcal{X}-\theta_\Phi\theta_A\|< 1 \quad \text{ and } \quad \|I_\mathcal{X}-\theta_\Psi\theta_B\|< 1.
	\end{align*} 
	In addition, if $( \{A_n\}_{n},  \{\Psi_n\}_{n} )$ and  $( \{B_n\}_{n},  \{\Phi_n\}_{n} )$ are   p-ASFs in $ \mathcal{B}(\mathcal{X}, \mathcal{Y})$, then we say that they are \textbf{approximately dual p-ASFs}.
\end{definition} 
 Definition \ref{APPROXIMATELYDUALPABS} gives the approximate duality using analysis and synthesis operators but  the following proposition  says that expansion of elements using approximate duals is close to the expansion using duals as described in Proposition  \ref{DUALCHAR}.
\begin{proposition}
	If $( \{A_n\}_{n},  \{\Psi_n\}_{n} )$ and  $( \{B_n\}_{n},  \{\Phi_n\}_{n} )$	are approximately dual operator-valued  p-ABSs in $ \mathcal{B}(\mathcal{X}, \mathcal{Y})$, then 
	\begin{align*}
		\left\|x-\sum_{n=1}^\infty
			\Phi_nA_nx \right\|< \|x\| \quad \text{ and } \quad \left\|x-\sum_{n=1}^\infty
		\Psi_nB_n x \right\|< \|x\|, \quad \forall x \in \mathcal{X}\setminus\{0\}.
	\end{align*}
\end{proposition}
\begin{proof}
	Given $ x \in \mathcal{X}\setminus\{0\}$, 
	\begin{align*}
		& \left\|x-\sum_{n=1}^\infty
		\Phi_nA_nx \right\|=\|I_\mathcal{X}x-\theta_\Phi \theta_A x\|\leq \|I_\mathcal{X}-\theta_\Phi \theta_A\|\|x\|< \|x\|, \\
		&	  \left\|x-\sum_{n=1}^\infty
			\Psi_nB_nx \right\|=\|I_\mathcal{X}x-\theta_\Psi\theta_Bx\|\leq \|I_\mathcal{X}-\theta_\Psi\theta_B\|\|x\|< \|x\|.
	\end{align*}
\end{proof}
One of the most important properties  associated with approximate duals is that they generate duals.
\begin{proposition}
	Let  $( \{A_n\}_{n},  \{\Psi_n\}_{n} )$ and  $( \{B_n\}_{n},  \{\Phi_n\}_{n} )$ be approximately dual  operator-valued  p-ABSs  	in $ \mathcal{B}(\mathcal{X}, \mathcal{Y})$. Then $ (\{B_nS_{B,\Psi}^{-1}\}_{n}, \{S_{A,\Phi}^{-1}\Phi_n\}_{n}) $ is a dual operator-valued  p-ABS for $( \{A_n\}_{n},  \{\Psi_n\}_{n} )$ 	and $ (\{A_nS_{A,\Phi}^{-1}\}_{n}, \{S_{B,\Psi}^{-1}\Psi_n\}_{n}) $ is a dual operator-valued  p-ABS for $( \{B_n\}_{n},  \{\Phi_n\}_{n} )$. 
\end{proposition}
\begin{proof}
Let $x \in \mathcal{X}$. Then 	
\begin{align*}
	&\sum_{n=1}^{\infty}\Psi_nB_nS_{B,\Psi}^{-1}x=S_{B,\Psi}S_{B,\Psi}^{-1}x=x, \quad 	\sum_{n=1}^{\infty}S_{A,\Phi}^{-1}\Phi_nA_nx=S_{A,\Phi}^{-1}S_{A,\Phi}x=x,\\
	&\sum_{n=1}^{\infty}\Phi_nA_nS_{A,\Phi}^{-1}x=S_{A,\Phi}S_{A,\Phi}^{-1}x=x, \quad \sum_{n=1}^{\infty}S_{B,\Psi}^{-1}\Psi_nB_nx=S_{B,\Psi}^{-1}S_{B,\Psi}x=x.
\end{align*}
\end{proof}
Motivated from  a characterization of duals, we try to get an operator-theoretic characterization of approximate duality.
\begin{theorem}\label{APPROXIMATEBESSELCHAR}
	Let $( \{A_n\}_{n},  \{\Psi_n\}_{n} )$ be an operator-valued  p-ABS in $ \mathcal{B}(\mathcal{X}, \mathcal{Y})$. 	An operator-valued p-ABS  $( \{B_n\}_{n},  \{\Phi_n\}_{n} )$ in $ \mathcal{B}(\mathcal{X}, \mathcal{Y})$ is an approximately dual operator-valued  p-ABS   for $( \{A_n\}_{n},  \{\Psi_n\}_{n} )$ if and only if there exist bounded linear operators  $U,V:	\mathcal{X} \to \mathcal{X}$ satisfying 	$\|I_\mathcal{X}-U\|< 1$ and $\|I_\mathcal{X}-V\|< 1$ such that $ (\{C_n\coloneqq B_nU^{-1}\}_{n}, \{\Xi_n\coloneqq V^{-1}\Phi_n\}_{n}) $ is a dual for  $( \{A_n\}_{n},  \{\Psi_n\}_{n} )$. Statement holds even if p-ABS is replaced by p-ASF.
\end{theorem}
\begin{proof}
	$(\Rightarrow)$ Define $U\coloneqq \theta_\Psi\theta_B$ and $V\coloneqq \theta_\Phi\theta_A$. We then have $\|I_\mathcal{X}-U\|< 1$ and   $\|I_\mathcal{X}-V\|< 1$ which say that $U$ and $V$ are invertible, and hence 
	\begin{align*}
		&\sum_{n=1}^\infty \Psi_nC_nx=\sum_{n=1}^\infty \Psi_nB_nU^{-1}x=\theta_\Psi \theta_BU^{-1}x =UU^{-1}x=x, \\
		&\sum_{n=1}^{\infty}\Xi_n A_nx=\sum_{n=1}^{\infty}V^{-1}\Phi_n A_nx=V^{-1}\theta_\Phi \theta_Ax=V^{-1}Vx=x, \quad \forall x \in
		\mathcal{X}.
	\end{align*}
	Hence $ (\{C_n\}_{n}, \{\Xi_n\}_{n}) $ is a dual for $( \{A_n\}_{n},  \{\Psi_n\}_{n} )$.
	
	$(\Leftarrow)$ We see  that $\theta_C=\theta_BU^{-1}$ and $\theta_\Xi=V^{-1}\theta_\Phi$. Since $ (\{C_n\}_{n}, \{\Xi_n\}_{n}) $ is a dual for $( \{A_n\}_{n},  \{\Psi_n\}_{n} )$, we get that 
	\begin{align*}
		&\|I_\mathcal{X}-\theta_\Phi\theta_A\|	=\|I_\mathcal{X}-V\theta_\Xi\theta_A\|=\|I_\mathcal{X}-V\|<1,\\
		&\|I_\mathcal{X}-\theta_\Psi\theta_B\|	=\|I_\mathcal{X}-\theta_\Psi\theta_CU\|=\|I_\mathcal{X}-U\|<1.
	\end{align*}
	Hence $( \{B_n\}_{n},  \{\Phi_n\}_{n} )$  is an approximately dual operator-valued p-ABS   for $( \{A_n\}_{n},  \{\Psi_n\}_{n} )$.
\end{proof}
Theorem  \ref{APPROXIMATEBESSELCHAR}   combined with Theorem  \ref{ALLDUALPOVF}   gives the following result. 
\begin{theorem}
	Let $( \{A_n\}_{n},  \{\Psi_n\}_{n} )$ be an operator-valued  p-ASF in $ \mathcal{B}(\mathcal{X}, \mathcal{Y})$. 	An operator-valued  p-ASF  $( \{B_n\}_{n},  \{\Phi_n\}_{n} )$ in $ \mathcal{B}(\mathcal{X}, \mathcal{Y})$ is an approximately dual operator-valued  p-ASF   for $( \{A_n\}_{n},  \{\Psi_n\}_{n} )$ if and only if there exist bounded linear operators  $U,V:	\mathcal{X} \to \mathcal{X}$, $A:\mathcal{X} \to \ell^p(\mathbb{N})$, $ B:\ell^p(\mathbb{N})\to \mathcal{X}$ satisfying 	$\|I_\mathcal{X}-U\|< 1$, $\|I_\mathcal{X}-V\|< 1$, 
\begin{align*}
	&B_n=A_nS_{A,\Psi}^{-1}U+\Gamma_nAU-A_nS_{A,\Psi}^{-1}\theta_\Psi AU,\\
	&\Phi_n=VS_{A,\Psi}^{-1}\Psi_n+VBL_n-VB\theta_AS_{A,\Psi}^{-1}\Psi_n, \quad \forall n \in \mathbb{N}
\end{align*}
such that the operator 
\begin{align*}
	S_{A, \Psi}^{-1}+BA-B\theta_AS_{A, \Psi}^{-1}\theta_\Psi A
\end{align*}
is bounded invertible, where   $A :\mathcal{X} \to \ell^p(\mathbb{N})\otimes\mathcal{Y}$ and $B:\ell^p(\mathbb{N})\otimes\mathcal{Y}\to\mathcal{X}$ are bounded linear operators.	
	\end{theorem}
In \cite{CHRISTENSENLAUGESEN},   Christensen and Laugesen  gave a method to construct approximately duals   iteratively.   In \cite{KRISHNAJOHNSON4}, this result was derived for Banach spaces. Here we have a  similar result for operator-valued p-ABSs.
\begin{theorem}
	Let 	$( \{A_n\}_{n},  \{\Psi_n\}_{n} )$ be an operator-valued  p-ABS in $ \mathcal{B}(\mathcal{X}, \mathcal{Y})$ and $( \{B_n\}_{n},  \{\Phi_n\}_{n} )$ be an operator-valued  p-ABS in $ \mathcal{B}(\mathcal{X}, \mathcal{Y})$ which is an approximately dual for $( \{A_n\}_{n},  \{\Psi_n\}_{n} )$.
	\begin{enumerate}[\upshape(i)]
		\item The dual operator-valued p-ABS  $ (\{B_nS_{B,\Psi}^{-1}\}_{n}, \{S_{A,\Phi}^{-1}\Phi_n\}_{n}) $   for $( \{A_n\}_{n},  \{\Psi_n\}_{n} )$ can be written as 
		\begin{align*}
			&B_nS_{B,\Psi}^{-1}=B_n+\sum_{m=1}^{\infty}B_n(I_\mathcal{X}-S_{B,\Psi})^m,\\
		&S_{A,\Phi}^{-1}\Phi_n= \Phi_n+\sum_{m=1}^{\infty}(I_\mathcal{X}-S_{A,\Phi})^m\Phi_n, \quad \forall n \in \mathbb{N}.
		\end{align*}
		\item For $N\in \mathbb{N}$, define 
		\begin{align*}
			&C_n^{(N)}\coloneqq  B_n+\sum_{m=1}^{N}B_n(I_\mathcal{X}-S_{B,\Psi})^m,\\
			&\Xi_n^{(N)}\coloneqq \Phi_n+\sum_{m=1}^{N}(I_\mathcal{X}-S_{A, \Phi})^m\Phi_n, \quad \forall n \in \mathbb{N}.
		\end{align*}
		Then $ (\{C_n^{(N)}\}_{n}, \{\Xi_n^{(N)}\}_{n}) $ is an operator-valued  p-ABS and is an approximately dual for $( \{A_n\}_{n},$ $  \{\Psi_n\}_{n} )$, for each $N\in \mathbb{N}$. Moreover, 
		\begin{align*}
			&\|I_\mathcal{X}-\theta_\Xi^{(N)}\theta_A\|\leq \|I_\mathcal{X}-S_{A,\Phi}\|^{N+1}\to 0 \quad \text{ as } N \to \infty, \\
			&	\|I_\mathcal{X}-\theta_\Psi\theta_C^{(N)}\|\leq \|I_\mathcal{X}-S_{B,\Psi}\|^{N+1}\to 0 \quad \text{ as } N \to \infty.
		\end{align*}
	\end{enumerate}
\end{theorem}
\begin{proof}
	\begin{enumerate}[\upshape(i)]
		\item We get these  from the Neumann series expansion 
		\begin{align*}
			S_{B,\Psi}^{-1}=\sum_{m=0}^{\infty}(I_\mathcal{X}-S_{B,\Psi})^m, \quad S_{A,\Phi}^{-1}=\sum_{m=0}^{\infty}(I_\mathcal{X}-S_{A,\Phi})^m.
		\end{align*}
		\item Clearly $ (\{C_n^{(N)}\}_{n}, \{\Xi_n^{(N)}\}_{n}) $ is an operator-valued  p-ABS. Now consider
		\begin{align*}
			\theta_\Xi^{(N)}\theta_Ax&=\sum_{n=1}^{\infty}\Xi_n^{(N)}A_nx=\sum_{n=1}^{\infty}\sum_{m=0}^{N}(I_\mathcal{X}-S_{A,\Phi})^m\Phi_nA_nx\\
			&=\sum_{m=0}^{N}\sum_{n=1}^{\infty}(I_\mathcal{X}-S_{A,\Phi})^m\Phi_nA_nx=\sum_{m=0}^{N}(I_\mathcal{X}-S_{A,\Phi})^m\sum_{n=1}^{\infty}\Phi_nA_nx\\
			&=\sum_{m=0}^{N}(I_\mathcal{X}-S_{A,\Phi})^mS_{A,\Phi}x=\sum_{m=0}^{N}(I_\mathcal{X}-S_{A,\Phi})^m(I_\mathcal{X}-(I_\mathcal{X}-S_{A,\Phi}))x\\
		&=x-(I_\mathcal{X}-S_{A,\Phi})^{N+1}x,\quad \forall x \in \mathcal{X}\\
			&\implies \|I_\mathcal{X}-\theta_\Xi^{(N)}\theta_A\|\leq \|I_\mathcal{X}-S_{A,\Phi}\|^{N+1}
		\end{align*}
		and 
		\begin{align*}
			\theta_\Psi\theta_C^{(N)}x&=	\sum_{n=1}^{\infty}\Psi_nC_n^{(N)}x=\sum_{n=1}^{\infty}\sum_{m=0}^N\Psi_nB_n(I_\mathcal{X}-S_{B,\Psi})^mx\\
			&=\sum_{m=0}^N\sum_{n=1}^\infty \Psi_nB_n(I_\mathcal{X}-S_{B,\Psi})^mx=\sum_{m=0}^NS_{B,\Psi}(I_\mathcal{X}-S_{B,\Psi})^mx\\
			&=\sum_{m=0}^N(I_\mathcal{X}-(I_\mathcal{X}-S_{B,\Psi}))(I_\mathcal{X}-S_{B,\Psi})^mx=x-(I_\mathcal{X}-S_{B,\Psi})^{N+1}x, \quad  \forall x \in \mathcal{X}\\
			&\implies \|I_\mathcal{X}-\theta_\Psi\theta_C^{(N)}\|\leq \|I_\mathcal{X}-S_{B,\Psi}\|^{N+1}.
		\end{align*}
		Since  $\|I_\mathcal{X}-S_{A,\Phi}\|<1$ and $\|I_\mathcal{X}-S_{B,\Psi}\|<1$, 	 $ (\{C_n^{(N)}\}_{n}, \{\Xi_n^{(N)}\}_{n}) $ is an approximately dual  for $( \{A_n\}_{n},  \{\Psi_n\}_{n} )$.
	\end{enumerate}
\end{proof}
In   \cite{CHRISTENSENLAUGESEN},  it is also showed by Christensen and Laugesen  that by perturbing frame we can get  approximate duals. In \cite{KRISHNAJOHNSON4}, it is showed that this result is valid for   Banach spaces. Here is the operator-valued version of that result.
\begin{theorem}
	Let $\{A_n\}_n$ be a collection in $ \mathcal{B}(\mathcal{X}, \mathcal{Y})$ and $\{\Psi_n\}_n$ be a collection in $ \mathcal{B}(\mathcal{Y}, \mathcal{X})$.  	Let  $ (\{C_n\}_{n}, \{\Xi_n\}_{n}) $ be an operator-valued  p-ASF in $ \mathcal{B}(\mathcal{X}, \mathcal{Y})$	  such that for some $R>0,$ $ Q>0$
	\begin{align}\label{ANALYSISEST}
		\left\|\sum_{n=1}^{\infty}L_n(A_n-C_n)x\right\|\leq R\|x\|,\quad \forall x \in \mathcal{X}
	\end{align}
	and 
	\begin{align}\label{SYNTHESISEST}
		\left\|\sum_{n=1}^{\infty}(\Xi_n-\Psi_n)\Gamma_nz\right\| 	\leq Q\|z\|,\quad  \forall z \in \ell^p(\mathbb{N})\otimes \mathcal{Y}.
	\end{align}
	Let $( \{B_n\}_{n},  \{\Phi_n\}_{n} )$ be a dual frame for $ (\{C_n\}_{n}, \{\Xi_n\}_{n}) $ with analysis bound $c$  and synthesis bound $d$. If $dR<1,$ $cQ<1$, 
	then  $( \{B_n\}_{n},  \{\Phi_n\}_{n} )$ is an approximately dual  for $ (\{A_n\}_{n}, \{\Psi_n\}_{n}) $.
\end{theorem}
\begin{proof}
	Inequalities (\ref{ANALYSISEST}) and (\ref{SYNTHESISEST}) say that $ (\{A_n\}_{n}, \{\Psi_n\}_{n}) $ is an operator-valued  p-ABS in $ \mathcal{B}(\mathcal{X}, \mathcal{Y})$.	Now  Inequality (\ref{ANALYSISEST}) can be written as $\|\theta_Ax-\theta_Cx\|\leq R\|x\|$, $\forall x \in \mathcal{X}$. Similarly, Inequality (\ref{SYNTHESISEST}) can be written as $\|\theta_\Xi z-\theta_\Psi z\|\leq Q\|z\|, \forall z  \in \ell^p(\mathbb{N})\otimes \mathcal{Y}$. Therefore 
	\begin{align*}
		\|I_\mathcal{X}-\theta_\Phi\theta_A\|=\|\theta_\Phi\theta_C-\theta_\Phi\theta_A\|\leq \|\theta_\Phi\|\|\theta_C-\theta_A\|\leq dR< 1, 
	\end{align*}	
	and 
	\begin{align*}
		\|I_\mathcal{X}-\theta_\Psi\theta_B\|=\|\theta_\Xi\theta_B-\theta_\Psi\theta_B\|\leq \|\theta_\Xi-\theta_\Psi\|\|\theta_B\|\leq cQ<1.
	\end{align*}
\end{proof}

\section{Equivalence}\label{SIMILARITYCOMPOSITIONANDTENSORPRODUCT}
Given two frames for Hilbert space, one tries to get one from another.  It is also clear that action of every element of frame by a fixed invertible bounded linear operator again gives a frame. This led Balan to introduce the notion of similarity or equivalence for frames for Hilbert spaces \cite{RADU}. This notion has been generalized in \cite{KAFTAL} for operator-valued frames, in \cite{KRISHNAJOHNSON5}  for factorable weak-OVFs  and in\cite{KRISHNAJOHNSON} for p-ASFs. Here we define the same notion for operator-valued p-ASFs.
\begin{definition}
	An   operator-valued p-ASF $( \{B_n\}_{n},  \{\Phi_n\}_{n} )$  in $\mathcal{B}(\mathcal{X}, \mathcal{Y})$    is said to be \textbf{similar}  or \textbf{equivalent} to an  operator-valued p-ASF  $  ( \{A_n\}_{n},  \{\Psi_n\}_{n} ) $ in $\mathcal{B}(\mathcal{X}, \mathcal{Y})$   if there exist   invertible operators  $ R_{A,B}, L_{\Psi, \Phi} $ $\in \mathcal{B}(\mathcal{X})$   such that 
	\begin{align*}
	B_n=A_nR_{A,B} , \quad \Phi_n=L_{\Psi, \Phi}\Psi_n, \quad \forall n \in \mathbb{N}.
	\end{align*} 
\end{definition}
We first observe that 	for every operator-valued p-ASF $  ( \{A_n\}_{n},  \{\Psi_n\}_{n} ) $, each  of operator-valued p-ASF $( \{A_nS_{A, \Psi}^{-1}\}_{n}, \{\Psi_n\}_{n})$    and  $ (\{A_n \}_{n}, \{S_{A,\Psi}^{-1}\Psi_n\}_{n})$ is a Parseval operator-valued p-ASF which is similar to  $  ( \{A_n\}_{n},  \{\Psi_n\}_{n} ) $.  
Like in the case of Hilbert spaces, since  $ R_{A,B}, R_{\Psi, \Phi}$ are bounded invertible,  the relation ``similarity" is an equivalence relation on the set 
\begin{align*}
\{( \{A_n\}_{n},  \{\Psi_n\}_{n} ):( \{A_n\}_{n},  \{\Psi_n\}_{n} ) \text{  is an operator-valued p-ASF in } \mathcal{B}(\mathcal{X}, \mathcal{Y})\}.
\end{align*} 
Next lemma shows that we can easily get analysis, synthesis and frame operators of similar frames from one another.
\begin{lemma}\label{SIM}
	Let $  ( \{A_n\}_{n},  \{\Psi_n\}_{n} ) $ and  $  ( \{B_n\}_{n},  \{\Phi_n\}_{n} ) $ be similar operator-valued p-ASFs in  $\mathcal{B}(\mathcal{X}, \mathcal{Y})$ and   $B_n=A_nR_{A,B} $, $\Phi_n=L_{\Psi, \Phi}\Psi_n,  \forall n \in \mathbb{N}$, for some invertible $ R_{A,B} ,L_{\Psi, \Phi} \in \mathcal{B}(\mathcal{X}).$ Then 
	\begin{enumerate}[\upshape(i)]
		\item $ \theta_B=\theta_A R_{A,B}, \theta_\Phi=L_{\Psi,\Phi}\theta_\Psi $.
		\item $S_{B,\Phi}=L_{\Psi,\Phi}S_{A, \Psi}R_{A,B}$.
		\item $P_{B,\Phi}=P_{A, \Psi}.$
	\end{enumerate}
\end{lemma}
\begin{proof}
	Let $x \in \mathcal{X}$ and $z \in \ell^p(\mathbb{N})\otimes \mathcal{Y}$. Then 
$ \theta_Bx=\sum_{n=1}^\infty L_nB_nx=\sum_{n=1}^\infty L_nA_nR_{A,B}x=\theta_AR_{A,B}x $ and  $ \theta_\Phi z=\sum_{n=1}^{\infty}\Phi_n\Gamma_n z=\sum_{n=1}^{\infty}L_{\Psi, \Phi}\Psi_n\Gamma_nz=L_{\Psi, \Phi}\theta_\Psi z$. Therefore 
\begin{align*}
	S_{B,\Phi}=\theta_\Phi\theta_B=L_{\Psi,\Phi}\theta_\Psi\theta_A R_{A,B}=L_{\Psi,\Phi}S_{A, \Psi}R_{A,B}
\end{align*}
and 
\begin{align*}
	P_{B, \Phi}=\theta_BS_{B,\Phi}^{-1}\theta_\Phi=(\theta_AR_{A,B})(L_{\Psi,\Phi}S_{A, \Psi}R_{A,B})^{-1}(L_{\Psi,\Phi}\theta_\Psi)=P_{A,\Psi}.
\end{align*}
\end{proof}
In the next result we characterize similarity using operators which do not depend on natural numbers. Further, we give a formula for operators which give similarity. 
\begin{theorem}\label{RIGHTSIMILARITY}
For two operator-valued p-ASFs  $  ( \{A_n\}_{n},  \{\Psi_n\}_{n} ) $ and  $  ( \{B_n\}_{n},  \{\Phi_n\}_{n} ) $ in $\mathcal{B}(\mathcal{X}, \mathcal{Y})$, the following are equivalent.
	\begin{enumerate}[\upshape(i)]
		\item $B_n=A_nR_{A,B} $,  $\Phi_n=L_{\Psi, \Phi}\Psi_n ,  \forall n \in \mathbb{N},$ for some invertible  $ R_{A,B} ,L_{\Psi, \Phi} \in \mathcal{B}(\mathcal{X}). $
		\item $\theta_B=\theta_AR_{A,B} $, $  \theta_\Phi=L_{\Psi, \Phi}\theta_\Psi $ for some invertible  $ R_{A,B} ,L_{\Psi, \Phi} \in \mathcal{B}(\mathcal{X}). $
		\item $P_{B,\Phi}=P_{A,\Psi}.$
	\end{enumerate}
	If one of the above conditions is satisfied, then  invertible operators in  $ \operatorname{(i)}$ and  $ \operatorname{(ii)}$ are unique and are given by $R_{A,B}=S_{A,\Psi}^{-1}\theta_\Psi\theta_B$,   $L_{\Psi, \Phi}=\theta_\Phi\theta_A S_{A,\Psi}^{-1}.$ In the case that $  ( \{A_n\}_{n},  \{\Psi_n\}_{n} ) $ is Parseval, then $  ( \{B_n\}_{n},  \{\Phi_n\}_{n} ) $ is  Parseval if and only if $R_{\Psi, \Phi}L_{A,B}=I_\mathcal{X} $   if and only if $L_{A,B}R_{\Psi, \Phi}=I_\mathcal{X} $.
\end{theorem}
\begin{proof}
 Lemma \ref{SIM} gives (i) $\Rightarrow$ (ii) $\Rightarrow$ (iii).  We show  (ii) $\Rightarrow$ (i). Assume (ii). We try to recover $B_n$ and $\Phi_n$. Using Equation (\ref{LEQUATION}), 
 \begin{align*}
 	&B_n x=\Gamma_n \theta_Bx=\Gamma_n \theta_AR_{A,B}x=A_n R_{A,B}x, \quad \forall x \in \mathcal{X}, \forall n \in \mathbb{N},\\
 	&\Phi_nz=\theta_\Phi L_n z=L_{\Psi, \Phi}\theta_\Psi L_nz=L_{\Psi, \Phi}\Psi_n , \quad  \forall n \in \mathbb{N}, \forall z \in \ell^p(\mathbb{N})\otimes \mathcal{Y}.
 \end{align*}
  We show (iii) $\Rightarrow$ (ii). Assume (iii). A calculation says that  $ \theta_B=P_{B,\Phi}\theta_B$ and $ \theta_\Phi=\theta_\Phi P_{B,\Phi}$. This observation says that 
  \begin{align*}
  	&\theta_B=P_{B,\Phi}\theta_B=P_{A,\Psi}\theta_B=\theta_A(S_{A,\Psi}^{-1}\theta_\Psi\theta_B)\\
  	&\theta_\Phi=\theta_\Phi P_{B,\Phi}=\theta_\Phi P_{A,\Psi}=(\theta_\Phi\theta_A S_{A,\Psi}^{-1})\theta_\Psi.
  \end{align*}
  To obtain (ii) it suffices to show that both $S_{A,\Psi}^{-1}\theta_\Psi\theta_B$ and $\theta_\Phi\theta_A S_{A,\Psi}^{-1}$ are invertible. For, 
  \begin{align*}
  	&(S_{A,\Psi}^{-1}\theta_\Psi\theta_B)(S_{B,\Phi}^{-1}\theta_\Phi\theta_A)=S_{A,\Psi}^{-1}\theta_\Psi P_{B,\Phi}\theta_A=S_{A,\Psi}^{-1}\theta_\Psi P_{A,\Psi}\theta_A=S_{A,\Psi}^{-1}\theta_\Psi \theta_A=I_\mathcal{X},\\
  	&(S_{B,\Phi}^{-1}\theta_\Phi\theta_A)(S_{A,\Psi}^{-1}\theta_\Psi\theta_B)=S_{B,\Phi}^{-1}\theta_\Phi P_{A,\Psi}\theta_B=S_{B,\Phi}^{-1}\theta_\Phi P_{B,\Phi}\theta_B=S_{B,\Phi}^{-1}\theta_\Phi\theta_B=I_\mathcal{X}
  \end{align*}
  and 
  \begin{align*}
  &(\theta_\Phi\theta_A S_{A,\Psi}^{-1})(\theta_\Psi\theta_BS_{B,\Phi}^{-1})=\theta_\Phi P_{A,\Psi}\theta_BS_{B,\Phi}^{-1}=\theta_\Phi P_{B,\Phi}\theta_BS_{B,\Phi}^{-1}=\theta_\Phi \theta_BS_{B,\Phi}^{-1}=I_\mathcal{X},\\
  &(\theta_\Psi\theta_BS_{B,\Phi}^{-1})(\theta_\Phi\theta_A S_{A,\Psi}^{-1})=\theta_\Psi P_{B,\Phi}\theta_AS_{A,\Psi}^{-1}=\theta_\Psi P_{A,\Psi}\theta_AS_{A,\Psi}^{-1}=\theta_\Psi \theta_AS_{A,\Psi}^{-1}=I_\mathcal{X}.
  \end{align*}
  	Let $ R_{A,B}, L_{\Psi,\Phi} \in \mathcal{B}(\mathcal{X}) $ be invertible. We saw that   $ R_{A,B}$ and $L_{\Psi,\Phi} $ satisfy (i) if and only if  they satisfy (ii). Let $B_n=A_nR_{A,B} $, $ \Phi_n=L_{\Psi, \Phi}\Psi_n ,  \forall n \in \mathbb{N}.$ Using (ii) $\theta_B=\theta_AR_{A,B} ,$ $ \theta_\Phi= L_{\Psi, \Phi}\theta_\Psi$. Therefore $\theta_\Psi\theta_B=\theta_\Psi\theta_AR_{A,B}=S_{A,\Psi}R_{A,B} $ and  $ \theta_\Phi \theta_A= L_{\Psi, \Phi}\theta_\Psi\theta_A=L_{\Psi, \Phi}S_{A,\Psi}$ which give the required formulas. 
\end{proof}
\begin{corollary}
	For any given operator-valued p-ASF $  ( \{A_n\}_{n},  \{\Psi_n\}_{n} ) $, the canonical dual of $  ( \{A_n\}_{n}, $ $ \{\Psi_n\}_{n} ) $ is the only dual operator-valued p-ASF  that is similar to $  ( \{A_n\}_{n},  \{\Psi_n\}_{n} ) $.
\end{corollary}
\begin{proof}
	Let  $ (\{B_n\}_{n} , \{\Phi_n\}_{n} )$ be an operator-valued p-ASF which is both dual and similar to  $  ( \{A_n\}_{n},  \{\Psi_n\}_{n} ) $.  Then  we have $ \theta_\Psi\theta_B=I_\mathcal{X}=\theta_\Phi\theta_A$ and  there exist invertible $ R_{A,B},L_{\Psi,\Phi}\in \mathcal{B}(\mathcal{X})$ such that  $B_n=A_nR_{A,B} ,$ $ \Phi_n=L_{\Psi, \Phi}\Psi_n ,$ $  \forall n \in \mathbb{N} $. Theorem \ref{RIGHTSIMILARITY} gives  $R_{A,B}=S_{A,\Psi}^{-1}\theta_\Psi\theta_B, $ $L_{\Psi, \Phi}=\theta_\Phi\theta_A S_{A,\Psi}^{-1}.$ But then $R_{A,B}=S_{A,\Psi}^{-1}I_\mathcal{X}=S_{A,\Psi}^{-1}$,  $ L_{\Psi, \Phi}=I_\mathcal{X}S_{A,\Psi}^{-1}=S_{A,\Psi}^{-1}.$ Therefore  $ (\{B_n\}_{n} , \{\Phi_n\}_{n} )$ is the  canonical  dual for  $  ( \{A_n\}_{n},  \{\Psi_n\}_{n} ) $.
\end{proof}
\begin{corollary}
	Two similar operator-valued p-ASFs   cannot be orthogonal.
\end{corollary}
\begin{proof}
	Let an operator-valued p-ASF  $  ( \{B_n\}_{n},  \{\Phi_n\}_{n} ) $ be similar to $  ( \{A_n\}_{n},  \{\Psi_n\}_{n} ) $. Then there are invertible  $ R_{A,B},L_{\Psi,\Phi}\in \mathcal{B}(\mathcal{X})$ such that  $B_n=A_nR_{A,B} , $ $\Phi_n=L_{\Psi, \Phi}\Psi_n ,  \forall n \in \mathbb{N} $. Using 	Theorem \ref{RIGHTSIMILARITY} and the invertibility  of  $R_{A,B} $ and $S_{A,\Psi} $, we have
	\begin{align*}
	\theta_\Psi\theta_B=\theta_\Psi\theta_AR_{A,B}=S_{A,\Psi}R_{A,B}\neq 0.
	\end{align*}
\end{proof}

\section{Perturbations}\label{PERTURBATIONS}
 Most useful  Paley-Wiener theorem \cite{ARSOVE} in Hilbert spaces says that sequences which are close to orthonormal bases are Riesz bases and in Banach spaces it states  that sequences which are close to Schauder bases are again Schauder bases. This motivated the perturbation of frames for Hilbert spaces. First result about stability of frames for Hilbert spaces was derived by Christensen in \cite{PALEY1} and improved in papers  \cite{PALEY2, PALEY3}. For Banach frames,  framings and Schauder frames,  stability results are derived in  \cite{CHENLIZHENG, ZHUWANG, CHRISTENSENHEIL, STOEVA}. We obtained stability of p-ASFs in \cite{KRISHNAJOHNSON2}. Here we  derive stability results of  operator-valued p-ASFs. First we recall a result which is an improvement of   result of  Hilding \cite{HILDING}.
 \begin{theorem}\cite{CASAZZAKALTON, PALEY3, VANEIJNDHOVEN}\label{cc1}
	Let $ \mathcal{X}, \mathcal{Y}$ be Banach spaces, $ U : \mathcal{X}\rightarrow \mathcal{Y}$ be a bounded invertible operator. If  a bounded linear  operator $ V : \mathcal{X}\rightarrow \mathcal{Y}$ is  such that there exist  $ \alpha, \beta \in \left [0, 1  \right )$ with 
	$$ \|Ux-Vx\|\leq\alpha\|Ux\|+\beta\|Vx\|,\quad \forall x \in  \mathcal{X},$$
	then $ V $ is bounded invertible and 
	$$ \frac{1-\alpha}{1+\beta}\|Ux\|\leq\|Vx\|\leq\frac{1+\alpha}{1-\beta} \|Ux\|, \quad\forall x \in  \mathcal{X};$$
	$$ \frac{1-\beta}{1+\alpha}\frac{1}{\|U\|}\|y\|\leq\|V^{-1}y\|\leq\frac{1+\beta}{1-\alpha} \|U^{-1}\|\|y\|, \quad\forall y \in  \mathcal{Y}.$$
\end{theorem}

\begin{theorem}\label{OURPERTURBATION}
	Let $  ( \{A_n\}_{n},  \{\Psi_n\}_{n} ) $ be an operator-valued p-ASF in $\mathcal{B}(\mathcal{X}, \mathcal{Y})$. Assume that  a collection $ \{\Psi_n \}_{n}	$ in $\mathcal{B}(\mathcal{Y}, \mathcal{X})$  is such that there exist $\alpha, \beta, \gamma \geq 0$ with  $ \max\{\alpha+\frac{\gamma}{\sqrt{a}}, \beta\}<1$ and 
\begin{align}\label{PEREQUATIONA}
	\left\|\sum_{n=1}^{m}(\Psi_n-\Phi_n)\Gamma_n z\right\|\leq \alpha\left\|\sum_{n=1}^{m}\Psi_n\Gamma_nz\right \|+\gamma \|z\|+\beta\left\|\sum_{n=1}^{m}\Phi_n\Gamma_nz\right \|,  \quad \forall z  \in \ell^p(\mathbb{N})\otimes\mathcal{Y}, 	 m=1, \dots.
\end{align}	
Then $  ( \{A_n\}_{n},  \{\Phi_n\}_{n} ) $ is an  operator-valued p-ASF in $\mathcal{B}(\mathcal{X}, \mathcal{Y})$ with frame bounds 
\begin{align*}
\frac{\|S_{A,\Psi}^{-1}\|(1+\beta)}{1-(\alpha+\gamma\|\theta_A S_{A,\Psi}^{-1}\|)}\quad	\text{ and } \quad \frac{\|\theta_A\|((1+\alpha)\|\theta_\Psi\|+\gamma)}{1-\beta}.
\end{align*}
\end{theorem}
\begin{proof}
For $ m=1, \dots,$ 
\begin{align*}
	\left\| \sum\limits_{n=1}^m\Phi_n\Gamma_nz\right\|&\leq \left\| \sum\limits_{n=1}^m(\Phi_n-\Psi_n)\Gamma_nz\right\|+\left\| \sum\limits_{n=1}^m\Psi_n\Gamma_nz\right\|\\
&\leq(1+\alpha)\left\| \sum\limits_{n=1}^m\Psi_n \Gamma_n z\right\|+\beta\left\| \sum\limits_{n=1}^m\Phi_n\Gamma_nz\right\|+\gamma\|z\|,\quad \forall z  \in \ell^p(\mathbb{N})\otimes\mathcal{Y}. 	
\end{align*}
Hence 
\begin{align*}
	\left\| \sum\limits_{n=1}^m\Phi_n\Gamma_nz\right\|\leq	\frac{1+\alpha}{1-\beta}\left\| \sum\limits_{n=1}^m\Psi_n\Gamma_nz\right\|+\frac{\gamma}{1-\beta}\|z\|,  \quad \forall  z  \in \ell^p(\mathbb{N})\otimes\mathcal{Y}, 	m=1, \dots.
\end{align*}
Therefore $\theta_\Phi$ is  a well-defined bounded linear operator with norm estimate
\begin{align*}
 \|\theta_\Phi\|\leq \frac{1+\alpha}{1-\beta}\|\theta_\Psi\|+\frac{\gamma}{1-\beta}.
\end{align*} 
Now Equation (\ref{PEREQUATIONA}) gives 
\begin{align*}
	\left\|\sum_{n=1}^{\infty}(\Psi_n-\Phi_n)\Gamma_n z\right\|\leq \alpha\left\|\sum_{n=1}^{\infty}\Psi_n\Gamma_nz\right \|+\gamma \|z\|+\beta\left\|\sum_{n=1}^{\infty}\Phi_n\Gamma_nz\right \|,  \quad \forall z  \in \ell^p(\mathbb{N})\otimes\mathcal{Y}. 	
\end{align*}
This is same as 
\begin{align}\label{PEREQUATIONB}
	\|\theta_\Psi z-\theta_\Phi z\|\leq \alpha \|	\theta_\Psi z\|+\gamma \|z\|+\beta 	\|\theta_\Phi z\|, \quad \forall z  \in \ell^p(\mathbb{N})\otimes\mathcal{Y}. 	
\end{align}
For every $x\in \mathcal{X}$, by taking $z=\theta_AS_{A,\Psi}^{-1}x$ in Equation (\ref{PEREQUATIONB}), we get 
\begin{align*}
	\|\theta_\Psi \theta_AS_{A,\Psi}^{-1}x-\theta_\Phi \theta_AS_{A,\Psi}^{-1}x\|\leq \alpha \|	\theta_\Psi \theta_AS_{A,\Psi}^{-1}x\|+\gamma \|\theta_AS_{A,\Psi}^{-1}x\|+\beta 	\|\theta_\Phi \theta_AS_{A,\Psi}^{-1}x\|, \quad \forall x \in \mathcal{X}.	
\end{align*}
That is, 
\begin{align*}
	\|x-\theta_\Phi \theta_AS_{A,\Psi}^{-1}x\|&\leq \alpha \|	x\|+\gamma \|\theta_AS_{A,\Psi}^{-1}x\|+\beta 	\|\theta_\Phi \theta_AS_{A,\Psi}^{-1}x\|\\
	&\leq (\alpha+\gamma\|\theta_AS_{A,\Psi}^{-1}\|)\|x\|+\beta 	\|\theta_\Phi \theta_AS_{A,\Psi}^{-1}x\|, \quad \forall x \in \mathcal{X}.	
\end{align*}
Since $ \max\{\alpha+\gamma\|\theta_A S_{A,\Psi}^{-1}\|, \beta\}<1$,  Theorem \ref{cc1} says  that  operator $S_{A, \Phi}S_{A,\Psi}^{-1}$ is invertible and 
\begin{align*}
	\|(S_{A,\Phi} S_{A,\Psi}^{-1})^{-1}\| \leq \frac{1+\beta}{1-(\alpha+\gamma\|\theta_A S_{A,\Psi}^{-1}\|)}.
\end{align*}
Hence the operator $S_{A,\Phi}=(S_{A,\Phi}S_{A,\Psi}^{-1})S_{A,\Psi}$ is invertible. Therefore $ (\{A_n \}_{n}, \{\Phi_n \}_{n}) $ is an operator-valued  p-ASF in $\mathcal{B}(\mathcal{X}, \mathcal{Y})$. To get  frame bounds we calculate:
\begin{align*}
	&\| S_{A,\Phi}^{-1}\|\leq\|S_{A,\Psi}^{-1}\|\| S_{A,\Psi}S_{A,\Phi}^{-1}\| \leq \frac{\|S_{A,\Psi}^{-1}\|(1+\beta)}{1-(\alpha+\gamma\|\theta_A S_{A,\Psi}^{-1}\|)}\quad\text{ and }\\
	&\|S_{A,\Phi}\|\leq \|\theta_\Phi\|\|\theta_A\|\leq \frac{\|\theta_A\|((1+\alpha)\|\theta_\Psi\|+\gamma)}{1-\beta}.
\end{align*}

\end{proof}
\begin{theorem}
		Let $  ( \{A_n\}_{n},  \{\Psi_n\}_{n} ) $ be an operator-valued p-ASF in $\mathcal{B}(\mathcal{X}, \mathcal{Y})$. Assume that a collection $\{B_n \}_{n} $ in $\mathcal{B}(\mathcal{X}, \mathcal{Y})$ and a collection $ \{\Phi_n \}_{n}	$ in  $\mathcal{B}(\mathcal{Y}, \mathcal{X})$ are such that there exist $r,s,t,\alpha, \beta, \gamma \geq 0$ with  $ \max\{ \beta,s\}<1$ and
	\begin{align*}
		&\left\|\sum_{n=1}^{m}L_n(A_n-B_n)x\right\|\leq r\left\|\sum_{n=1}^{m}L_nA_nx\right \|+t \|x\|+s\left\|\sum_{n=1}^{m}L_nB_nx\right \|,   \quad\forall x  \in \mathcal{X}, m=1, \dots,\\
		&	\left\|\sum_{n=1}^{m}(\Psi_n-\Phi_n)\Gamma_n z\right\|\leq \alpha\left\|\sum_{n=1}^{m}\Psi_n\Gamma_nz\right \|+\gamma \|z\|+\beta\left\|\sum_{n=1}^{m}\Phi_n\Gamma_nz\right \|,  \quad \forall z  \in \ell^p(\mathbb{N})\otimes\mathcal{Y}, 	 m=1, \dots.
	\end{align*}	
	Assume that one of the following holds. 
	\begin{enumerate}
		\item $\sum_{n=1}^{\infty}\|(\Psi_nA_n-\Phi_nB_n)S_{A,\Psi}^{-1}\|<1.$
		\item $\sum_{n=1}^{\infty}\|S_{A,\Psi}^{-1}(\Psi_nA_n-\Phi_nB_n)\|<1.$
		\item $\sum_{n=1}^{\infty}\|S_{A,\Psi}^{-1}\Psi_nA_n-\Phi_nB_nS_{A,\Psi}^{-1}\|<1.$
		\item $\sum_{n=1}^{\infty}\|\Psi_nA_nS_{A,\Psi}^{-1}-S_{A,\Psi}^{-1}\Phi_nB_n\|<1$.
	\end{enumerate}
	Then $ (\{B_n \}_{n}, \{\Phi_n \}_{n}) $ is an operator-valued p-ASF in $\mathcal{B}(\mathcal{X}, \mathcal{Y})$. Moreover, an upper frame  bound is  
	\begin{align*}
		\left(\frac{1+\alpha}{1-\beta}\|\theta_\Psi\|+\frac{\gamma}{1-\beta}\right)\left(\frac{1+r}{1-s}\|\theta_A\|+\frac{t}{1-s}\right).
	\end{align*}
\end{theorem}
\begin{proof}
	Similar to the arguments given in the beginning of  proof of Theorem \ref{OURPERTURBATION}, we see that $\theta_B$ and $\theta_\Phi$ are well-defined bounded linear operators. 	We deal with  four cases.\\
	Assume (1). Then 
	\begin{align*}
		\left\|x-\sum_{n=1}^{\infty}\Phi_nB_nS_{A,\Psi}^{-1}x\right\|&=\left\|\sum_{n=1}^{\infty}\Psi_nA_nS_{A,\Psi}^{-1}x-\sum_{n=1}^{\infty}\Phi_nB_nS_{A,\Psi}^{-1}x\right\|\\
		&=\left\|\sum_{n=1}^{\infty}(\Psi_nA_nS_{A,\Psi}^{-1}-\Phi_nB_nS_{A,\Psi}^{-1})x\right\|\\
		&\leq \sum_{n=1}^{\infty}\|(\Psi_nA_n-\Phi_nB_n)S_{A,\Psi}^{-1}\|\|x\|, \quad \forall x \in \mathcal{X}.
	\end{align*}
	Therefore the operator  $S_{B,\Phi}S_{A,\Psi}^{-1}$ is invertible.\\
	Assume (2). Then 
	\begin{align*}
		\left\|x-\sum_{n=1}^{\infty}S_{A,\Psi}^{-1}\Phi_nB_nx\right\|&=\left\|\sum_{n=1}^{\infty}S_{A,\Psi}^{-1}\Psi_nA_nx-\sum_{n=1}^{\infty}S_{A,\Psi}^{-1}\Phi_nB_nx\right\|\\
		&=\left\|\sum_{n=1}^{\infty}(S_{A,\Psi}^{-1}\Psi_nA_n-S_{A,\Psi}^{-1}\Phi_nB_n)x\right\|\\
		&\leq \sum_{n=1}^{\infty}\|S_{A,\Psi}^{-1}(\Psi_nA_n-\Phi_nB_n)\|\|x\|, \quad \forall x \in \mathcal{X}.
	\end{align*}
	Therefore the operator  $S_{A,\Psi}^{-1}S_{B,\Phi}$ is invertible.\\
	Assume (3). Then 
	\begin{align*}
	\left\|x-\sum_{n=1}^{\infty}\Phi_nB_nS_{A,\Psi}^{-1}x\right\|&=\left\|\sum_{n=1}^{\infty}S_{A,\Psi}^{-1}\Psi_nA_nx-\sum_{n=1}^{\infty}\Phi_nB_nS_{A,\Psi}^{-1}x\right\|\\
	&=\left\|\sum_{n=1}^{\infty}(S_{A,\Psi}^{-1}\Psi_nA_n-\Phi_nB_nS_{A,\Psi}^{-1})x\right\|\\
	&\leq \sum_{n=1}^{\infty}\|S_{A,\Psi}^{-1}\Psi_nA_n-\Phi_nB_nS_{A,\Psi}^{-1}\|\|x\|, \quad \forall x \in \mathcal{X}.
	\end{align*}
	Therefore the operator  $S_{B,\Phi}S_{A,\Psi}^{-1}$ is invertible.\\
	Assume (4). Then 
	\begin{align*}
		\left\|x-\sum_{n=1}^{\infty}S_{A,\Psi}^{-1}\Phi_nB_nx\right\|&=\left\|\sum_{n=1}^{\infty}\Psi_nA_nS_{A,\Psi}^{-1}x-\sum_{n=1}^{\infty}S_{A,\Psi}^{-1}\Phi_nB_nx\right\|\\
	&=\left\|\sum_{n=1}^{\infty}(\Psi_nA_nS_{A,\Psi}^{-1}-S_{A,\Psi}^{-1}\Phi_nB_n)x\right\|\\
	&\leq \sum_{n=1}^{\infty}\|\Psi_nA_nS_{A,\Psi}^{-1}-S_{A,\Psi}^{-1}\Phi_nB_n\|\|x\|, \quad \forall x \in \mathcal{X}.
	\end{align*}
	Therefore the operator  $S_{A,\Psi}^{-1}S_{B,\Phi}$ is invertible.
\end{proof} 
\section{Operator-valued Feichtinger conjectures}\label{FEICHTINGER}
Recently first author formulated  a variety of conjectures  and problems  for p-approximate Schauder frames for Banach spaces \cite{MAHESHKRISHNA, MAHESHKRISHNA2}. Here we present operator-valued versions of some of them.
\begin{conjecture}\label{FB}\textbf{(Feichtinger conjecture for operator-valued p-ASFs)
		Let  $ (\{A_n \}_{n},  \{\Psi_n \}_{n}) $  be  an operator-valued  p-ASF in  $\mathcal{B}(\mathcal{X}, \mathcal{Y})$ such that 
		\begin{align*}
			0<\inf_{n\in \mathbb{N}}\|\Psi_n\|\leq \sup_{n\in \mathbb{N}}\|\Psi_n\|<\infty \quad \text{ and } \quad 	0<\inf_{n\in \mathbb{N}}\|A_n\|\leq \sup_{n\in \mathbb{N}}\|A_n\|<\infty.
		\end{align*}
		Then  $ (\{A_n \}_{n}, \{\Psi_n \}_{n}) $  can be partitioned into a finite union of operator-valued  p-approximate Riesz sequences. Moreover, what is the number of partitions required?}
\end{conjecture}
\begin{conjecture}\label{FS}\textbf{(Feichtinger conjecture for operator-valued p-ABSs)
		Let   $ (\{A_n \}_{n}, \{\Psi_n \}_{n}) $  be  an operator-valued  p-ABS in  $\mathcal{B}(\mathcal{X}, \mathcal{Y})$   such that 
		\begin{align*}
			0<\inf_{n\in \mathbb{N}}\|\Psi_n\|\leq \sup_{n\in \mathbb{N}}\|\Psi_n\|<\infty \quad \text{ and } \quad 	0<\inf_{n\in \mathbb{N}}\|A_n\|\leq \sup_{n\in \mathbb{N}}\|A_n\|<\infty.
		\end{align*}
		Then $ (\{A_n \}_{n}, \{\Psi_n \}_{n}) $  can be partitioned into a finite union of operator-valued  p-approximate Riesz sequences. Moreover, what is the number of partitions required?}
\end{conjecture}

\begin{definition}
	A 	operator-valued p-ASF $ (\{A_n \}_{n}, \{\Psi_n \}_{n}) $  in $\mathcal{B}(\mathcal{X}, \mathcal{Y})$   is said to be \textbf{p-scalable} if there exist  sequences of  operators $\{T_n\}_n$ and $\{S_n\}_n$   in $\mathcal{B}(\mathcal{X})$   such that 
	\begin{align*}
		 (\{T_nC_n \}_{n}, \{S_n\Psi_n \}_{n}) 	 \text{  is a Parseval  operator-valued  p-ASF in  } \mathcal{B}(\mathcal{X}, \mathcal{Y}).
	\end{align*}	
\end{definition}
\begin{problem}
	(\textbf{Scaling problem for Banach spaces}) \textbf{For given Banach spaces $\mathcal{X}$,  $\mathcal{Y}$, classify operator-valued p-ASFs $ (\{A_n \}_{n}, \{\Psi_n \}_{n}) $  in   $ \mathcal{B}(\mathcal{X}, \mathcal{Y})$  so that there are  sequences of  operators  $\{T_n\}_n$, $\{S_n\}_n$ in $\mathcal{B}(\mathcal{X})$   such that $(\{A_nT_n\}_n, \{S_n\Psi_n\}_n)$ is a Parseval operator-valued p-ASF   in $ \mathcal{B}(\mathcal{X}, \mathcal{Y})$, i.e.,  $ (\{A_n \}_{n}, \{\Psi_n \}_{n}) $ is p-scalable.}
\end{problem}
  \textbf{Note that, due to Theorem 3.1 in  \cite{SUN1} (and due to the fundamental  works of Marcus, Spielman, Srivastava, Casazza, Edidin,  Weaver, Fickus, Tremain, Vershynin, Anderson, and Weber   \cite{MARCUSSPIELMANSRIVASTAVA, CASAZZAEDIDIN, WEAVER, CASAZZAFICKUSTREMAINWEBER, CASAZZATREMAINKADISONSINGER, ANDERSON}), operator-valued Feichtinger conjectures for operator-valued frames/G-frames are true for Hilbert spaces. }
  \section{Acknowledgements}
 First author thanks  Prof. B. V. Rajarama Bhat for the Post Doc position through his J. C. Bose Fellowship (SERB).

 \bibliographystyle{plain}
 \bibliography{reference.bib}

\end{document}